\documentclass[10pt]{article}
\usepackage{amssymb,amsmath,amsthm}
\usepackage[numbers,sort&compress]{natbib}
\usepackage{indentfirst}
\usepackage{exscale}
\usepackage{relsize}
\usepackage{subfigure}
\usepackage{graphicx}

\newcommand{\R}{\mathbb{R}}
\theoremstyle{definition}

\theoremstyle{remark}

\numberwithin{equation}{section}
\allowdisplaybreaks
 \UseRawInputEncoding

\begin{document}
\title{\Large\bf{ Existence of three solutions for a poly-Laplacian system on graphs}}
\date{}
\author {Yan Pang$^1$, \ Xingyong Zhang$^{1,2}$\footnote{Corresponding author, E-mail address: zhangxingyong1@163.com}\\
{\footnotesize $^1$Faculty of Science, Kunming University of Science and Technology,}\\
 {\footnotesize Kunming, Yunnan, 650500, P.R. China.}\\
{\footnotesize $^{2}$Research Center for Mathematics and Interdisciplinary Sciences, Kunming University of Science and Technology,}\\
 {\footnotesize Kunming, Yunnan, 650500, P.R. China.}\\}

 \date{}
 \maketitle

 \begin{center}
 \begin{minipage}{15cm}
 \par
 \small  {\bf Abstract:} We deal with the existence of three distinct solutions for a poly-Laplacian system with a parameter on  finite  graphs  and a $(p,q)$-Laplacian system  with a parameter  on locally finite graphs.  The main tool is an abstract critical points theorem  in [Bonanno and Bisci, J.Math.Appl.Anal, 2011, 382(1): 1-8]. A key point in this paper is that we overcome the difficulty to prove that the G$\hat{\mbox{a}}$teaux derivative of the variational functional for poly-Laplacian operator admits a continuous inverse, which is caused by  the special definition of the poly-Laplacian operator on graph and mutual coupling of two variables in system.

 \end{minipage}
 \end{center}
  \allowdisplaybreaks
 \vskip2mm
 {\section{Introduction }}
\setcounter{equation}{0}

Let $G=(V, E)$ be a graph with the vertex set  $V$  and   the edge set $E$. If both $V$ and $E$ are finite set, then $G$ is called as a finite graph. If for any $x\in V$, there are finite  vertexes $y\in V$ such that $xy\in E$ ($xy$ denotes an edge connecting $x$ with $y$), then $G$ is called as a locally finite graph.  For any edge $xy\in E$ with two vertexes $x,y$, let $\omega_{xy}(>0)$ denote its weight and suppose that $\omega_{xy}=\omega_{yx}$.  For any $x\in V$, let $deg(x)=\sum_{y\thicksim x}\omega_{xy}$, where $y\thicksim x$ denotes those $y$ connecting $x$ with $xy\in E$. Suppose that $\mu:V\rightarrow \R^+$ is a finite measure.
Define the directional derivative of a function $u:V\rightarrow\mathbb{R}$ by
\begin{eqnarray}
\label{Eq5}
D_{w,y}u(x):=\frac{1}{\sqrt2}(u(x)-u(y))\sqrt\frac{w_{xy}}{\mu(x)}.
\end{eqnarray}
Define the gradient of $u$  as a vector
\begin{eqnarray}
\label{Eq6}
\nabla u(x):=(D_{w,y}u(x))_{y\in V},
\end{eqnarray}
which is indexed by the vertices $y\in V$.
It is easy to obtain that $\nabla(u_1+u_2)=\nabla u_1+\nabla u_2$ and
\begin{eqnarray*}
\label{Eq9}
\nabla u\nabla v
&  =  & (D_{w,y}u(x))_{y\in V}(D_{w,y}v(x))_{y\in V}\\
&  =  &  \sum\limits_{y\thicksim x}D_{w,y}u(x)D_{w,y}v(x)\\
&  =  &  \frac{1}{2\mu(x)}\sum\limits_{y\thicksim x}w_{xy}(u(y)-u(x))(v(y)-v(x)).
\end{eqnarray*}
Let
\begin{eqnarray}
\label{Eq7}
\Gamma(u,v)(x)=\frac{1}{2\mu(x)}\sum\limits_{y\thicksim x}w_{xy}(u(y)-u(x))(v(y)-v(x)).
\end{eqnarray}
Then
\begin{eqnarray}
\label{Eq8}
\Gamma(u,v)=\nabla u\nabla v.
\end{eqnarray}
Let $\Gamma(u)=\Gamma(u,u)$ and define the length of the gradient by
\begin{eqnarray}
\label{Eq10}
|\nabla u|(x)=\sqrt{\Gamma(u)(x)}=\left(\frac{1}{2\mu(x)}\sum\limits_{y\thicksim x}w_{xy}(u(y)-u(x))^2\right)^{\frac{1}{2}}.
\end{eqnarray}
The Laplacian operator of $u:V\to\R$ is defined by
\begin{eqnarray}
\label{Eq4}
\Delta u(x):=\frac{1}{\mu(x)}\sum\limits_{y\thicksim x}w_{xy}(u(y)-u(x)).
\end{eqnarray}
Let $|\nabla^m u|$ denote the length of $m$-order gradient of $u$, which is defined by
\begin{eqnarray}
\label{Eq11}
|\nabla^mu|=
 \begin{cases}
  |\nabla\Delta^{\frac{m-1}{2}}u|,& \text {when $m$ is odd,}\\
  |\Delta^{\frac{m}{2}}u|,&  \text {when $m$ is even,}
   \end{cases}
\end{eqnarray}
where $\nabla\Delta^{\frac{m-1}{2}}u$ is defined as in $(\ref{Eq6})$ with $u$ replaced by $\Delta^{\frac{m-1}{2}}u$, and  $\Delta^{\frac{m}{2}}u$ is defined by  $\Delta^{\frac{m}{2}}u=\Delta(\Delta^{\frac{m}{2}-1}u)$ which means that $u$ is replaced by $\Delta^{\frac{m}{2}-1}u$ in (\ref{Eq4}).
 By mathematical induction, we can obtain that
 \begin{eqnarray}\label{aa1}
 \Delta^{\frac{m}{2}}(u_1+u_2)=\Delta^{\frac{m}{2}}u_1+\Delta^{\frac{m}{2}} u_2,\ \ \mbox{if } m \mbox{ is even.}
 \end{eqnarray}

\par
For any given $p>1$, the $p$-Laplacian operator is defined by
\begin{eqnarray}\label{Eq13}
\Delta_pu(x):=\frac{1}{2\mu(x)}\sum\limits_{y\sim x}\left(|\nabla u|^{p-2}(y)+|\nabla u|^{p-2}(x)\right)\omega_{xy}(u(y)-u(x)).
\end{eqnarray}
It is easy to see that $p$-Laplacian operator reduces to the Laplacian operator of $u$ if $p=1$.
\par
For any function $u:V\rightarrow\mathbb{R}$, we denote
\begin{eqnarray}\label{Eq12}
\int\limits_V u(x) d\mu=\sum\limits_{x\in V}\mu(x)u(x).
\end{eqnarray}
For any given real number $r\ge 1$, we define
$$
L^r(V)=\left\{u:V\to\R\Big|\int_V|u(x)|^rd\mu<\infty\right\}
$$
endowed with the norm
\begin{eqnarray}
\label{Eq2}
\|u\|_{L^r(V)}=\left(\int_V|u(x)|^rd\mu\right)^\frac{1}{r}.
\end{eqnarray}
In the distributional sense, $\Delta_p u$ can be written as follows. For any $u\in\mathcal{C}_c(V)$,
\begin{eqnarray}\label{Eq14}
\int\limits_V(\Delta_p u)v d\mu=-\int\limits_V|\nabla u|^{p-2}\Gamma(u,v)d\mu,
\end{eqnarray}
where $\mathcal{C}_c(V)$ is the set of all real functions with compact support.
Furthermore,  a more general operator can be introduced, denoted by $\pounds_{m,p}$,  as following: for any function $\phi:V\to\R$,
\begin{eqnarray}\label{eq9}
\int\limits_V(\pounds_{m,p}u)\phi d\mu=
 \begin{cases}
  \int_V|\nabla^m u|^{p-2}\Gamma(\Delta^{\frac{m-1}{2}}u,\Delta^{\frac{m-1}{2}}\phi)d\mu,& \text { if $m$ is odd},\\
  \int_V|\nabla^m u|^{p-2}\Delta^{\frac{m}{2}}u\Delta^{\frac{m}{2}}\phi d\mu,&  \text { if $m$ is even}.
    \end{cases}
\end{eqnarray}
where $m\ge 1$ are integers and $p>1$.  The operator $\pounds_{m,p}$ is called as the poly-Laplacian operator of $u$ if $p=2$, and obviously, the operator $\pounds_{m,p}$ reduces to the $p$-Laplacian operator if $m=1$. The above knowledge are taken from   \cite{Chung2005} and \cite{Yamabe 2016}.

\par
 In this paper, we study
 the existence of three  solutions for the following  poly-Laplacian system
\begin{eqnarray}
\label{eq2}
 \begin{cases}
  \pounds_{m_1,p}u+h_1(x)|u|^{p-2}u=\lambda F_u(x,u,v),\;\;\;\;\hfill x\in V,\\
  \pounds_{m_2,q}v+h_2(x)|v|^{q-2}v=\lambda F_v(x,u,v),\;\;\;\;\hfill x\in V,\\
   \end{cases}
\end{eqnarray}
where $G=(V,E)$ is a finite graph, $m_i\geq1$ are integers, $h_i:V\rightarrow\mathbb{R},i=1,2$, $p,q>1$, $\lambda>0$, $F:V\times\mathbb{R}^2\rightarrow\mathbb{R}$ and $\pounds_{m_1,p}$ and $\pounds_{m_2,q}$ are  defined by (\ref{eq9}).
\par
Moreover, if $G=(V,E)$ is a locally finite graph, we consider the existence of three  solutions for the following $(p,q)$-Laplacian system:
\begin{eqnarray}
\label{p1}
 \begin{cases}
   -\Delta_p u+h_1(x)|u|^{p-2}u=\lambda F_u(x,u,v),\;\;\;\;\hfill x\in V,\\
   -\Delta_q v+h_2(x)|v|^{q-2}v=\lambda F_v(x,u,v),\;\;\;\;\hfill x\in V,\\
 \end{cases}
\end{eqnarray}
where $-\Delta_p$ and $-\Delta_q$ are defined by (\ref{Eq13}) with $p\ge2$ and $q\ge2$, $F:V\times \R^2 \to \R$, $h_i:V\rightarrow\mathbb{R},i=1,2$,  and $\lambda> 0$.
\par
In recent years, some problems with respect to equations on  graphs attract some attentions. We refer readers to [5-7,10-12] as some examples.
In \cite{Yamabe 2016}, Grigior'yan-Lin-Yang investigated the following poly-Laplacian equation on  graph $G=(V,E)$:
\begin{eqnarray}
\label{(1)}
\pounds_{m,p}u+h(x)|u|^{p-2}u=\lambda f(x,u)\;\;\mbox{in}\;\;x\in V,
\end{eqnarray}
 where $p>1,h:V\rightarrow\mathbb{R}$ and $f:V\times\mathbb{R}\rightarrow\mathbb{R}$. They considered the case that the graph $G=(V,E)$ is a locally finite graph, $h(x)\equiv 0$ and the equation (\ref{(1)}) satisfies the Diriclet boundary condition, and the case that the graph $G=(V,E)$ is a finite graph. They established some existence results of a positive solution for equation (\ref{(1)}) with $\lambda=1$  by mountain pass theorem.
  In \cite{Grigor 2017}, Grigior'yan-Lin-Yang studied (\ref{(1)}) with $m=1$ and $p=2$, where $V$ is a locally finite graph. They obtained two existence results of positive solutions for equation (\ref{(1)}) by  mountain pass theorem.
 In \cite{Imbesi 2023}, when $m=1$ and $p=2$, by applying a three critical point theorem from \cite{Ricceri2000},  Imbesi-Bisci-Repovs established some existence results of at least two solutions for equation (\ref{(1)}) when the parameter $\lambda$ locates at some concrete range.
In \cite{Yu 2023}, Yu-Zhang-Xie-Zhang studied system (\ref{eq2}) with $\lambda=1$, $p=q$ and $F(x,u,v)$  satisfying asymptotically-$p$-linear growth at infinity with respect to $(u,v)$. By using the mountain pass theorem, they obtained that system (\ref{eq2}) has a nontrivial solution and they also presented some corresponding results for equation (\ref{(1)}) with $\lambda=1$.
 In \cite{Yang 2023}, Yang and Zhang investigated system (\ref{p1}) with perturbations and two parameters $\lambda_1$ and $\lambda_2$. When $F$ possesses sub-$(p, q)$ growth on $(u,v)$,  an existence result of one
nontrivial solution was established by Ekeland's variational principle, and when $F$  possesses  super-$(p, q)$ growth on $(u,v)$,  one solution of positive energy and one solution of negative energy were obtained by mountain pass theorem and Ekeland's variational principle, respectively.
  In \cite{Zhang 2022},
Zhang-Zhang-Xie-Yu considered system (\ref{eq2}) with $\lambda=1$. They established an existence result and a  multiplicity result of nontrivial solutions when $F$ satisfies the super-$(p,q)$ growth conditions on $(u,v)$ via  mountain pass theorem and  symmetric mountain pass theorem, respectively.

\par
In present paper, our work are mainly motivated by  \cite{Bonanno G2010}, \cite{Bonanno G2011} and \cite{Zhang 2022}. In \cite{Bonanno G2010}, Bonanno and Marano established an existence result of three critical points for $f_\lambda:=\Phi-\lambda\Psi$ with  $\lambda\in\mathbb{R}$, and got a
well-determined large interval of parameters  for which $f_\lambda$ possesses at least three
critical points under weaker regularity and compactness conditions. In \cite{Bonanno G2011}, Bonanno and Bisci obtained that a non-autonomous elliptic Dirichlet problem possesses at least three weak solutions by the three critical points theorem in \cite{Bonanno G2010}.
 \par
 Based on the works in \cite{Bonanno G2010}, \cite{Bonanno G2011} and \cite{Zhang 2022}, the motivation of our work is to consider whether the three critical points theorem due to Bonanno and Bisci in \cite{Bonanno G2010} can be applied to systems (\ref{eq2}) and (\ref{p1}). The main difficulty of such problem is to prove that the G$\hat{\mbox{a}}$teaux derivative of the variational functional $\Phi$ admits a continuous inverse. The main reason to cause this difficulty is that the special definition of $\pounds_{m,p}$ and mutual coupling of $u$ and $v$
in systems (\ref{eq2}) and (\ref{p1}) make proving the uniformly monotone of  $\Phi'$ become difficult. To overcome this difficulty, we discuss the case that $m$ is even and the case that $m$ is odd, respectively,  and sufficiently  use the formulars (\ref{Eq8}) and (\ref{aa1}), and in order to deal with mutual coupling of $u$ and $v$, we divide into four cases about the norms of $u$ and $v$ and then make some careful arguments (see Lemma 3.5 below).

We call that $(u,v)$ is a non-trivial solution of system   (\ref{eq2}) (or (\ref{p1})) if $(u,v)$ satisfies   (\ref{eq2}) (or (\ref{p1})) and $(u,v)\not=(0,0)$. Next, we state our results.

\vskip2mm
\noindent
{\bf (I) For the poly-Laplacian system on finite graph }
\vskip2mm
\noindent
{\bf Theorem 1.1.} {\it Let $G=(V,E)$ be a finite graph. Assume that the following conditions hold: \\
$(H)$\; $h_i(x)>0$ for all $x\in V$, $i=1,2$;\\
$(F_0)$\; $F(x,s,t)$ is continuously differentiable in $(s,t)\in \R^2$ for all $x\in V$;\\
$(F_1)$ \;   $\int_VF(x,0,0)d\mu=0$;\\
$(F_2)$ \; there exist two constants $\alpha \in [0,p),\beta \in [0,q)$ and functions $f_i,g:V\to \R, i=1,2$,  such that
$$
F(x,s,t)\le f_1(x)|s|^\alpha+ f_2(x)|t|^\beta+g(x)
$$
 for all $(x,s,t)\in V\times\mathbb{R}\times\mathbb{R}$;\\
$(F_3)$ \; there are  positive constants $\gamma_1$,  $\gamma_2$, $\delta_1$ and $\delta_2$ with $\delta_i>\gamma_i\kappa_i$, $i=1,2$, such that
$$
 \Lambda_1:=\frac{1}{\gamma_1^{p}+\gamma_2^{q}} \max_{x\in V,|s|\le{\frac{(p\gamma_1^{p}+p\gamma_2^{q})^{\frac{1}{p}}}{h_{1,\min}^{1/p}\mu_{\min}^{1/p}} }, |t|\le \frac{ (q\gamma_1^{p}+q\gamma_2^{q})^{\frac{1}{q}}}{h_{2,\min}^{1/q}\mu_{\min}^{1/q}}}F(x,s,t)|V| \\
< \frac{\inf_{x\in V}F(x,\delta_1,\delta_2)|V|}{ \frac{\delta_1^p}{p}\int_Vh_1(x)d\mu+\frac{\delta_2^q}{q}\int_Vh_2(x)d\mu}:=\Lambda_2,\\
$$
where $|V|=\sum_{x\in V}\mu(x)$, $h_{i,\min}=\min_{x\in V} h_i(x)$, $i=1,2$, $\mu_{\min}=\min_{x\in V}\mu(x)$,
$$
\kappa_1= \left(\frac{1}{p}\int_Vh_1(x)d\mu\right)^{-\frac{1}{p}}, \ \ \kappa_2=\left(\frac{1}{q}\int_Vh_2(x)d\mu\right)^{-\frac{1}{q}}.
$$
  Then  for each parameter $\lambda$ belonging to
$\left( \Lambda_2^{-1},  \Lambda_1^{-1} \right),$
 system (\ref{eq2}) has at least three distinct solutions.}

\vskip2mm
\noindent
{\bf (II) For the $(p,q)$-Laplacian system on locally finite graph }
\vskip2mm
\noindent
{\bf Theorem 1.2.} {\it  Let $G=(V,E)$ be a locally finite graph. Assume that the following conditions hold: \\
$(M)$\: there exists a $\mu_0>0$ such that $\mu(x)\ge \mu_0$ for all $x\in V$;\\
$(H_1)$\; there exists a constant $h_0>0$ such that $h_i(x)\geq h_0>0$ for all $x\in V$, $i=1,2$;\\
$(F_0)'$\; $F(x,s,t)$ is continuously differentiable in $(s,t)\in \R^2$ for all $x\in V$,  and there exists a  function $a\in C(\R^+,\R^+)$ and a function $b:V\to \R^+$ with $b\in L^1(V)$ such that
$$
|F_s(x,s,t)|\le a(|(s,t)|) b(x), |F_t(x,s,t)|\le a(|(s,t)|) b(x), |F(x,s,t)|\le a(|(s,t)|) b(x),
$$
for all $x\in V$ and all $(s,t)\in \R^2$;\\
 $(F_1)'$\; $\int_V F(x,0,0)d\mu=0$ and there exists a $x_0\in V$ such that $F(x_0,0,0)=0$;\\
$(F_2)'$ \; there exist two constants $\alpha \in [0,p),\beta \in [0,q)$ and functions $f_i,g:V\to \R, i=1,2$, with $f_i\in L^\infty(V)$, $i=1,2$ and $g\in L^1(V)$,  such that
$$
F(x,s,t)\leq f_1(x)|s|^\alpha+ f_2(x)|t|^\beta+g(x)
$$
 for all $(x,s,t)\in V\times\mathbb{R}\times\mathbb{R}$;\\
$(F_3)'$ \; there are  positive constants $\gamma_1$,  $\gamma_2$, $\delta_1$ and $\delta_2$ with $\delta_i>\gamma_i\kappa_i$, $i=1,2$, such that
\begin{eqnarray*}
&   & \Theta_1:=\frac{1}{\gamma_1^{p}+\gamma_2^{q}}\max_{|(s,t)|\le \frac{1}{h_0^{1/p}\mu_0^{1/p}} (p\gamma_1^{p}+p\gamma_2^{q})^{\frac{1}{p}}+ \frac{1}{h_0^{1/q}\mu_0^{1/q}} (q\gamma_1^{p}+q\gamma_2^{q})^{\frac{1}{q}}}a(|(s,t)|)\int_V b(x)d\mu \\
& < & \frac{F(x_0,\delta_1,\delta_2)}{ \frac{\delta_1^pM_1}{p}+\frac{\delta_2^q M_2}{q}}:=\Theta_2,
\end{eqnarray*}
where $\kappa_1= \left(\frac{M_1}{p}\right)^{-\frac{1}{p}}, \ \ \kappa_2=\left(\frac{M_2}{q}\right)^{-\frac{1}{q}}$, and
\begin{eqnarray*}
M_1& = & \left(\frac{deg(x_0)}{2\mu(x_0)}\right)^{\frac{p}{2}}\mu(x_0)+h_1(x_0)\mu(x_0)+\sum_{y\thicksim x_0}\left(\frac{w_{x_0y}}{2\mu(y)}\right)^{\frac{p}{2}}\mu(y),\\
M_2& = & \left(\frac{deg(x_0)}{2\mu(x_0)}\right)^{\frac{q}{2}}\mu(x_0)+h_2(x_0)\mu(x_0)+\sum_{y\thicksim x_0}\left(\frac{w_{x_0y}}{2\mu(y)}\right)^{\frac{q}{2}}\mu(y).
\end{eqnarray*}
  Then  for each parameter $\lambda$ belonging to
$\left( \Theta_2^{-1},  \Theta_1^{-1} \right)$,
 system (\ref{p1}) has at least three distinct solutions.}
\par

\vskip2mm
\noindent
{\bf Remark 1.1.} In Theorem 1.1 and Theorem 1.2, all of  three solutions are
nontrivial solutions if we furthermore assume that $F_s(x,0,0)\not=0$ or $F_t(x,0,0)\not=0$ for some $x\in V$.

\vskip2mm
{\section{Preliminaries}}
\setcounter{equation}{0}
\par
In this section, we recall the Sobolev space on graph and some embedding relationships (see \cite{Yamabe 2016,Zhang 2022, Yang 2023}). We also recall
an abstract critical point theorem in \cite{Bonanno G2010}, which is the main tool to prove our main results.
\par
Let $G=(V,E)$ be a graph.
For any given integer $m\geq1$ and any given real number $l>1$, we define
\begin{eqnarray*}
W^{m,l}(V)=\left\{u:V\to\R\Big|\int_V(|\nabla^m u(x)|^l+h(x)|u(x)|^l)d\mu<\infty\right\}
\end{eqnarray*}
 endowed with the norm
\begin{eqnarray*}
\label{Eq1}
\|u\|_{W^{m,l}(V)}=\left(\int_V(|\nabla^m u(x)|^l+h(x)|u(x)|^l)d\mu\right)^\frac{1}{l},
\end{eqnarray*}
where $h(x)>0$ for all $x\in V$. If $V$ is a finite graph, then $W^{m,l}(V)$ is of finite dimension.
\par

  \vskip2mm
\noindent
{\bf Lemma 2.1.} (\cite{Yamabe 2016, Zhang 2022}) {\it Let $G=(V,E)$ be a finite graph. For all $\psi\in W^{m,l}(V)$, there is
$$\|\psi\|_{\infty}\leq d_l\|\psi\|_{W^{m,l}(V)},$$
where $\|\psi\|_{\infty}=\max_{x\in V}|\psi(x)|$ and $d_l=\left(\frac{1}{\mu_{\min}h_{\min}}\right)^\frac{1}{l}$.}

 \vskip2mm
\noindent
{\bf Lemma 2.2.} (\cite{Yamabe 2016,Zhang 2022}) {\it Let $G=(V,E)$ be a finite graph. Then $W^{m,l}(V)\hookrightarrow L^r(V)$ for all $1\leq r\leq+\infty$. Especially, if $1< r<+\infty$, then for all $\psi\in W^{m,l}(V)$,
\begin{eqnarray*}\label{Eq15}
\|\psi\|_{L^r(V)}\leq K_{l,r}\|\psi\|_{W^{m,l}(V)},
\end{eqnarray*}
where}
$$
K_{l,r}=\frac{\left(\sum_{x\in V}\mu(x)\right)^{\frac{1}{r}}}{\mu_{\min}^{\frac{1}{l}}h_{\min}^{\frac{1}{l}}}.
$$

 \vskip2mm
\noindent
{\bf Lemma 2.3.} (\cite{Yang 2023}) {\it Let $G=(V,E)$ be a locally finite graph.  If $\mu(x)\geq\mu_0$ and $(H_1)$ holds, then $W_h^{1,l}(V)$ is continuously embedded into $L^r(V)$ for all $1<l\leq r\leq \infty$, and the following inequalities hold:
\begin{eqnarray*}\label{Eq16}
\|u\|_\infty\leq \frac{1}{h_0^{1/l}\mu_0^{1/l}}\|u\|_{W^{1,l}(V)}
\end{eqnarray*}
and}
\begin{eqnarray*}\label{Eq17}
\|u\|_{L^r(V)}\leq \mu_0^\frac{l-r}{lr}h_0^{-\frac{1}{l}}\|u\|_{W^{1,l}(V)}\;\;\mbox{for all}\; l\leq r<\infty.
\end{eqnarray*}

 \noindent
{\bf Lemma 2.4.} (\cite{Bonanno G2010}) {\it Let $W$ be a real reflexive Banach space, $\Phi:W\to\mathbb{R}$ be a coercive, continuously G$\hat{a}$teaux differentiable and sequentially weakly lower semicontinuous functional whose G$\hat{a}$teaux derivative admits a continuous inverse on $W^*$, $\Psi:W\to\mathbb{R}$ be a continuously G$\hat{a}$teaux differentiable functional whose G$\hat{a}$teaux derivative is compact such that
$$
\Phi(0)=\Psi(0)=0.
$$
Assume that there exist $r>0$ and $\bar{x}\in X$, with $r<\Phi(\bar{x})$, such that:\\
$(a_1)$ $\frac{\sup_{\Phi(x)\le{r}}\Psi(x)}{r}<\frac{\Psi(\bar{x})}{\Phi(\bar{x})}$;\\
$(a_2)$ for each $\lambda\in\Lambda_r:=\left(\frac{\Phi(\bar{x})}{\Psi(\bar{x})},\frac{r}{\sup_{\Phi(x)\le{r}}\Psi(x)}\right)$, the functional $\Phi-\lambda\Psi$ is coercive. \\
Then, for each $\lambda\in\Lambda_r$, the functional $\varphi_{\lambda}:=\Phi-\lambda\Psi$ has at least three distinct critical points in $W$.}

\vskip2mm
{\section{Proofs for the poly-Laplacian system (\ref{eq2})}}
  \setcounter{equation}{0}

  \par
 Let $G=(V,E)$ be a finite graph. In order to investigate the poly-Laplacian system (\ref{eq2}), we work in  the space $W:=W^{m_1,p}(V)\times W^{m_2,q}(V)$ with the norm  endowed with $\|(u,v)\|=\|u\|_{W^{m_1,p}(V)}+\|v\|_{W^{m_2,q}(V)}$. Then $(W,\|\cdot\|)$ is a  Banach space of finite dimension. Consider the functional $\varphi:W\to\R$ as
\begin{eqnarray}
\label{EQ1} \varphi(u,v)=\frac{1}{p}\int_V(|\nabla^{m_1}u|^p+h_1(x)|u|^p)d\mu+\frac{1}{q}\int_V(|\nabla^{m_2}v|^q+h_2(x)|v|^q)d\mu-\lambda\int_V F(x,u,v)d\mu.
\end{eqnarray}
Then under the assumptions of Theorem 1.1, $\varphi\in C^1(W,\R)$ and
\begin{eqnarray}\label{EQ2}
\langle\varphi'(u,v),(\phi_1,\phi_2)\rangle&=&\int_V\left[\pounds_{m_1,p}u\phi_1+h_1(x)|u|^{p-2}u\phi_1-\lambda  F_u(x,u,v)\phi_1\right]d\mu\nonumber\\
&&+\int_V\left[\pounds_{m_2,q}v\phi_2+h_2(x)|v|^{q-2}v\phi_2-\lambda F_v(x,u,v)\phi_2\right]d\mu
\end{eqnarray}
for any $(u,v),(\phi_1,\phi_2)\in W$.  In order to apply Lemma 2.4, we will use the functionals $\Phi:W\rightarrow\mathbb{R}$ and $\Psi:W\rightarrow\mathbb{R}$ defined by setting
  \begin{eqnarray*}
\Phi(u,v)&=&\frac{1}{p}\int_V(|\nabla^{m_1}u|^p+h_1(x)|u|^p)d\mu+\frac{1}{q}\int_V(|\nabla^{m_2}v|^q+h_2(x)|v|^q)d\mu\\
&=&\frac{1}{p}\|u\|^p_{W^{m_1,p}(V)}+\frac{1}{q}\|v\|^q_{W^{m_2,q}(V)}
  \end{eqnarray*}
  and
\begin{eqnarray*}
\Psi(u,v)=\int_V F(x,u,v)d\mu.
\end{eqnarray*}
Then $\varphi(u,v)=\Phi-\lambda\Psi$. Moreover, it is easy to see that $(u,v)\in W$ is a critical point of $\varphi$ if and only if
\begin{eqnarray*}
\int_V\left(\pounds_{m_1,p}u+h_1(x)|u|^{p-2}u-\lambda F_u(x,u,v)\right)\phi_1d\mu=0
\end{eqnarray*}
 and
\begin{eqnarray*}
\int_V\left(\pounds_{m_2,q}v+h_2(x)|v|^{q-2}v-\lambda F_v(x,u,v)\right)\phi_2d\mu=0
\end{eqnarray*}
for all $(\phi_1,\phi_2)\in W$. By the arbitraries of $\phi_1$ and $\phi_2$, we conclude that
\begin{eqnarray*}
\pounds_{m_1,p}u+h_1(x)|u|^{p-2}u=\lambda F_u(x,u,v),\\
\pounds_{m_2,q}v+h_2(x)|v|^{q-2}v=\lambda F_v(x,u,v).
\end{eqnarray*}
  Thus the problem of finding the solutions for system (\ref{eq2}) is reduced to seek the critical points of functional $\varphi$ on $W$.
  \par

 \vskip2mm
 \noindent
{\bf Lemma 3.1.}  {\it Assume that $(F_0)$ holds. Then for any given $r>0$, the following inequality holds:
$$
\frac{\sup_{(u,v)\in\Phi^{-1}(-\infty,r]}\Psi(u,v)}{r}\le \frac{1}{r}\max_{x\in V,|s|\le{\frac{(pr)^{\frac{1}{p}}}{h_{\min}^{1/p}\mu_{\min}^{1/p}} }, |t|\le \frac{(qr)^{\frac{1}{q}}}{h_{\min}^{1/q}\mu_{\min}^{1/q}} }F(x,u(x),v(x))|V|.
$$ }

 \vskip0mm
 \noindent
{\bf Proof.} By $(F_0)$, we have
\begin{eqnarray*}
         \Psi(u,v)&=&\int_V{F(x,u(x),v(x))}d\mu\\
&  =   & \sum_V F(x,u(x),v(x))d\mu\\
& \leq &  \max_{x\in V,|s|\le\|u\|_\infty, |t|\le \|v\|_\infty}F(x,s,t)|V|
\end{eqnarray*}
for every $(u,v)\in W$. Furthermore,  for all $(u,v)\in W$ with $\Phi(u,v)\leq r$, by Lemma 2.1, we get
\begin{eqnarray*}
\label{a1}
\|u\|_\infty \le \frac{1}{h_{1,\min}^{1/p}\mu_{\min}^{1/p}} \|u\|_{W^{m_1,p}(V)}\le \frac{1}{h_{1,\min}^{1/p}\mu_{\min}^{1/p}} (pr)^{\frac{1}{p}},\ \  \|v\|_\infty \le \frac{1}{h_{2,\min}^{1/q}\mu_{\min}^{1/q}} \|v\|_{W^{m_2,q}(V)}\le \frac{1}{h_{2,\min}^{1/q}\mu_{\min}^{1/q}} (qr)^{\frac{1}{q}}.
\end{eqnarray*}
Hence, we get
\begin{eqnarray}
\label{(q8)}
\sup_{(u,v)\in\Phi^{-1}(-\infty,r]}\Psi(u,v)\leq    \max_{x\in V,|s|\le{\frac{(pr)^{\frac{1}{p}}}{h_{1,\min}^{1/p}\mu_{\min}^{1/p}} }, |t|\le \frac{(qr)^{\frac{1}{q}}}{h_{2,\min}^{1/q}\mu_{\min}^{1/q}} }F(x,s,t)|V|.
\end{eqnarray}
Then the proof is completed by multiplying $\frac{1}{r}$ in both side of (\ref{(q8)}). \qed

\vskip2mm
 \noindent
{\bf Lemma 3.2.} {\it Assume that $(F_0)$ and $(F_3)$ hold. Then there exists $ (u_{\delta_1},v_{\delta_2})\in W$  such that}
$$
 \frac{\sup_{(u,v)\in\Phi^{-1}(-\infty,\gamma_1^{p}+\gamma_2^{q}]}\Psi(u,v)}{\gamma_1^{p}+\gamma_2^{q}}
 < \frac{\Psi(u_{\delta_1},v_{\delta_2})}{\Phi(u_{\delta_1},v_{\delta_2})}.
$$
\vskip0mm
 \noindent
{\bf Proof.}   Let
\begin{eqnarray*}
u_{\delta_{1}}(x)=\delta_1,
v_{\delta_{2}}(x)=\delta_{2},\ \ \forall x\in V,
\end{eqnarray*}
where $\delta_i, i=1,2$ are given in $(F_3)$. It is easy to verify that  $(u_{\delta_1},v_{\delta_2})\in W$, $|\nabla^{m_1}u_{\delta_1}|=0$ and $|\nabla^{m_1}v_{\delta_1}|=0$ for all $m_i\ge 1,i=1,2.$
Then
\begin{eqnarray}
\label{(q5)}
    \Phi(u_{\delta_1},v_{\delta_2})
 & = &  \frac{1}{p}\int_V{(|\nabla^{m_1}u_{\delta_1}|^p+h_1(x)|u_{\delta_1}|^p)}d\mu+\frac{1}{q}\int_V{(|\nabla^{m_2}v_{\delta_2}|^q+h_2(x)|v_{\delta_2}|^q)}d\mu\nonumber\\
 & = &  \frac{\delta_1^p}{p}\int_Vh_1(x)d\mu+\frac{\delta_2^q}{q}\int_Vh_2(x)d\mu.
\end{eqnarray}
Note that $\delta_i>\gamma_i\kappa_i$, $i=1,2$, where
$$
\kappa_1= \left(\frac{1}{p}\int_Vh_1(x)d\mu\right)^{-\frac{1}{p}}, \ \ \kappa_2=\left(\frac{1}{q}\int_Vh_2(x)d\mu\right)^{-\frac{1}{q}}.
$$
Then (\ref{(q5)}) implies that $ \Phi(u_{\delta_1},v_{\delta_2})\ge \gamma_1^p+\gamma_2^q$.
 Moreover,
\begin{eqnarray}
\label{(q6)}
        \Psi(u_{\delta_1},v_{\delta_2})&=&\int_VF(x,u_{\delta_1},v_{\delta_2})d\mu\nonumber\\
&  =   & \int_{V}F(x,\delta_1,\delta_2)d\mu\nonumber\\
& \geq & \inf_{x\in V}F(x,\delta_1,\delta_2)|V|.
\end{eqnarray}
Hence, by (\ref{(q5)}) and (\ref{(q6)}), we get
\begin{eqnarray}\label{qq1}
      \frac{\Psi(u_{\delta_1},v_{\delta_2})}{\Phi(u_{\delta_1},v_{\delta_2})}
\geq  \frac{\inf_{x\in V}F(x,\delta_1,\delta_2)|V|}{ \frac{\delta_1^p}{p}\int_Vh_1(x)d\mu+\frac{\delta_2^q}{q}\int_Vh_2(x)d\mu}.
\end{eqnarray}
In view of Lemma 3.1, (\ref{qq1}) and ($F_3$), we get
\begin{eqnarray*}
&      &  \frac{\sup_{(u,v)\in\Phi^{-1}(-\infty,\gamma_1^{p}+\gamma_2^{q}]}\Psi(u,v)}{\gamma_1^{p}+\gamma_2^{q}}\\
& \leq & \frac{1}{\gamma_1^{p}+\gamma_2^{q}} \max_{x\in V,|s|\le{\frac{1}{h_{\min}^{1/p}\mu_{\min}^{1/p}} (p(\gamma_1^{p}+\gamma_2^{q}))^{\frac{1}{p}}}, |t|\le \frac{1}{h_{\min}^{1/q}\mu_{\min}^{1/q}} (q(\gamma_1^{p}+\gamma_2^{q}))^{\frac{1}{q}}}F(x,s,t)|V| \\
& < & \frac{\inf_{x\in V}F(x,\delta_1,\delta_2)|V|}{ \frac{\delta_1^p}{p}\int_Vh_1(x)d\mu+\frac{\delta_2^q}{q}\int_Vh_2(x)d\mu}\\
& \leq & \frac{\Psi(u_{\delta_1},v_{\delta_2})}{\Phi(u_{\delta_1},v_{\delta_2})}.
\end{eqnarray*}
The proof is completed. \qed

\vskip2mm
\noindent
{\bf Lemma 3.3.} {\it Assume that $(F_2)$ holds. Then for each $\lambda\in \left(0,+\infty\right)$, the functional $\Phi-\lambda\Psi$ is coercive.}
\vskip0mm
\noindent
{\bf Proof.} By $(F_2)$, we have
\begin{eqnarray*}
       \varphi(u,v)
&  =   & \frac{1}{p}\|u\|_{W^{m_1,p}(V)}^p+\frac{1}{q}\|v\|_{W^{m_2,q}(V)}^q-\lambda\int_V{F(x,u,v)}d\mu\\
& \ge  & \frac{1}{p}\|u\|_{W^{m_1,p}(V)}^p+\frac{1}{q}\|v\|_{W^{m_2,q}(V)}^q-\lambda \|f_1\|_\infty\|u\|_{\infty}^{\alpha}-\lambda \|f_2\|_\infty\|v\|_{\infty}^\beta-\lambda\int_Vg(x)d\mu\\
& \ge  & \frac{1}{p}\|u\|_{W^{m_1,p}(V)}^p+\frac{1}{q}\|v\|_{W^{m_2,q}(V)}^q-\lambda \|f_1\|_\infty d_p^\alpha\|u\|_{W^{m_1,p}(V)}^\alpha\\
&      &  -\lambda \|f_2\|_\infty d_q^\beta\|v\|_{W^{m_2,q}(V)}^\beta -\lambda\int_Vg(x)d\mu.
\end{eqnarray*}
Note that $\alpha\in[0,p)$ and $\beta\in[0,q)$. Therefore $\varphi(u,v)$ is a coercive functional for every $\lambda\in \left(0,+\infty\right)$.
\qed

\vskip2mm
\noindent
{\bf Lemma 3.4.} {\it The G$\hat{a}$teaux derivative of $\Phi$ admits a continuous inverse on $W^*$, where  $W^*$ is the dual space of $W$.}
\vskip0mm
\noindent
{\bf Proof.} Firstly, we prove that $\Phi'$ is uniformly monotone. Via $(2.2)$ of \cite{Simon J 1977}, there exists a positive constant $c_p$ such that
\begin{eqnarray}
\label{(11)}
(|x|^{p-2}x-|y|^{p-2}y,x-y)\geq c_p|x-y|^p,\ \ \mbox{for all } x,y\in \R^N,
\end{eqnarray}
where $(\cdot,\cdot)$ denotes the inner product in $\R^N$.
Note that
\begin{eqnarray*}
&&\langle\Phi'(u_1,v_1)-\Phi'(u_2,v_2),(u_1-u_2,v_1-v_2)\rangle\\
&  =   & \int_V\left[\pounds_{m_1,p}u_1(u_1-u_2)+h_1(x)|u_1|^{p-2}u_1(u_1-u_2)\right]d\mu\\
&  +   & \int_V\left[\pounds_{m_2,q}v_1(v_1-v_2)+h_2(x)|v_1|^{q-2}v_1(v_1-v_2)\right]d\mu\\
&  -   & \int_V\left[\pounds_{m_1,p}u_2(u_1-u_2)+h_1(x)|u_2|^{p-2}u_2(u_1-u_2)\right]d\mu\\
&  -   & \int_V\left[\pounds_{m_2,q}v_2(v_1-v_2)+h_2(x)|v_2|^{q-2}v_2(v_1-v_2)\right]d\mu.
\end{eqnarray*}
Firstly, we prove that
\begin{eqnarray*}
   \uppercase\expandafter{\romannumeral1}
& : =  &\int_V\left[\pounds_{m_1,p}u_1(u_1-u_2)+h_1(x)|u_1|^{p-2}u_1(u_1-u_2)\right.\\
&     &\left.-\pounds_{m_1,p}u_2(u_1-u_2)+h_1(x)|u_2|^{p-2}u_2(u_1-u_2)\right]d\mu\\
&\geq &c_p\|u_1-u_2\|^p_{W^{m_1.p}(V)}.
\end{eqnarray*}
When $m_1$ is odd, by (\ref{eq9}), (\ref{Eq8}), (\ref{Eq11}), (\ref{aa1}) and (\ref{(11)}), we have
\begin{eqnarray*}
&&\int_V\left[\pounds_{m_1,p}u_1(u_1-u_2)-\pounds_{m_1,p}u_2(u_1-u_2)\right]d\mu\\
&  =   &  \int_V\left[|\nabla^{m_1} u_1|^{p-2}\Gamma(\Delta^\frac{m_1-1}{2}u_1,\Delta^\frac{m_1-1}{2}(u_1-u_2))-
|\nabla^{m_1} u_2|^{p-2}\Gamma(\Delta^\frac{m_1-1}{2}u_2,\Delta^\frac{m_1-1}{2}(u_1-u_2))\right]d\mu\\
&  =   &  \int_V\left[|\nabla^{m_1} u_1|^{p-2}\nabla\Delta^\frac{m_1-1}{2}u_1\cdot\nabla\Delta^\frac{m_1-1}{2}(u_1-u_2)-
|\nabla^{m_1} u_2|^{p-2}\nabla\Delta^\frac{m_1-1}{2}u_2\cdot\nabla\Delta^\frac{m_1-1}{2}(u_1-u_2)\right]d\mu\\
&  =   &  \int_V\left[|\nabla\Delta^\frac{m_1-1}{2}u_1| ^{p-2}\nabla\Delta^\frac{m_1-1}{2}u_1\cdot\nabla\Delta^\frac{m_1-1}{2}(u_1-u_2)-
|\nabla\Delta^\frac{m_1-1}{2}u_2|^{p-2}\nabla\Delta^\frac{m_1-1}{2}u_2\cdot\nabla\Delta^\frac{m_1-1}{2}(u_1-u_2)\right]d\mu\\
&  =   &  \int_V\left[\nabla\Delta^\frac{m_1-1}{2}(u_1-u_2)\cdot(|\nabla\Delta^\frac{m_1-1}{2}u_1| ^{p-2}\nabla\Delta^\frac{m_1-1}{2}u_1-|\nabla\Delta^\frac{m_1-1}{2}u_2|^{p-2}\nabla\Delta^\frac{m_1-1}{2}u_2)\right]d\mu\\
& \geq &  \int_V c_p|\nabla\Delta^\frac{m_1-1}{2}(u_1-u_2)|^pd\mu=  \int_Vc_p|\nabla^{m_1}(u_1-u_2)|^pd\mu.
\end{eqnarray*}
When $m_1$ is even, by (\ref{eq9}), (\ref{Eq11}), (\ref{aa1}) and (\ref{(11)}), we also have
\begin{eqnarray*}
&&\int_V\left[\pounds_{m_1,p}u_1(u_1-u_2)-\pounds_{m_1,p}u_2(u_1-u_2)\right]d\mu\\
&  =   &  \int_V\left[|\nabla^{m_1} u_1|^{p-2}\Delta^\frac{m_1}{2}u_1\Delta^\frac{m_1}{2}(u_1-u_2)-
|\nabla^{m_1} u_2|^{p-2}\Delta^\frac{m_1}{2}u_2\Delta^\frac{m_1}{2}(u_1-u_2)\right]d\mu\\
&  =   &  \int_V\left[|\Delta^\frac{m_1}{2} u_1|^{p-2}\Delta^\frac{m_1}{2}u_1\Delta^\frac{m_1}{2}(u_1-u_2)-
|\Delta^\frac{m_1}{2} u_2|^{p-2}\Delta^\frac{m_1}{2}u_2\Delta^\frac{m_1}{2}(u_1-u_2)\right]d\mu\\
&  =   &  \int_V\left[\Delta^\frac{m_1}{2}(u_1-u_2)(|\Delta^\frac{m_1}{2} u_1|^{p-2}\Delta^\frac{m_1}{2}u_1-
|\Delta^\frac{m_1}{2} u_2|^{p-2}\Delta^\frac{m_1}{2}u_2)\right]d\mu\\
& \geq &  \int_Vc_p|\Delta^\frac{m_1}{2} (u_1-u_2)|^pd\mu= \int_Vc_p|\nabla^{m_1}(u_1-u_2)|^pd\mu.
\end{eqnarray*}
Thus, for all positive integer $m$, we get
\begin{eqnarray}
\label{(9)}
\int_V\left[\pounds_{m_1,p}u_1(u_1-u_2)-\pounds_{m_1,p}u_2(u_1-u_2)\right]d\mu\geq
 \int_Vc_p|\nabla^{m_1}(u_1-u_2)|^pd\mu.
\end{eqnarray}
By (\ref{(11)}), we have
\begin{eqnarray}
\label{(10)}
\int_V\left[h_1(x)|u_1|^{p-2}u_1(u_1-u_2)-h_1(x)|u_2|^{p-2}u_2(u_1-u_2)\right]d\mu\geq \int_Vh_1(x)c_p|u_1-u_2|^pd\mu.
\end{eqnarray}
So, by (\ref{(9)}) and (\ref{(10)}), we obtain that
\begin{eqnarray*}
\uppercase\expandafter{\romannumeral1}\geq c_p\|u_1-u_2\|^p_{W^{m_1.p}(V)}.
\end{eqnarray*}
Similarly, we can prove that there exists a positive constant $c_q$ such that
\begin{eqnarray*}
\uppercase\expandafter{\romannumeral2}
& : =  &\int_V\left[\pounds_{m_2,q}v_1(v_1-v_2)+h_2(x)|v_1|^{q-2}v_1(v_1-v_2)\right.\\
&     &\left.-\pounds_{m_2,q}v_2(v_1-v_2)+h_2(x)|v_2|^{q-2}v_2(v_1-v_2)\right]d\mu\\
&\geq &c_q\|v_1-v_2\|^q_{W^{m_2.q}(V)}.
\end{eqnarray*}
Hence,
\begin{eqnarray}
\label{eqq4}
\langle\Phi'(u_1,v_1)-\Phi'(u_2,v_2),(u_1-u_2,v_1-v_2)\rangle\nonumber\\
\geq c_p\|u_1-u_2\|^p_{W^{m_1.p}(V)}
+c_q\|v_1-v_2\|^q_{W^{m_2.q}(V)}.
\end{eqnarray}
Next, we consider the following four cases if we let $\max\{p,q\}=p$.
\par
(1) Assume that $\|u_1-u_2\|_{W^{m_1.p}(V)}>1$  and $\|v_1-v_2\|_{W^{m_2.q}(V)}>1$. Then $\|(u_1-u_2,v_1-v_2)\|>2$ and
\begin{eqnarray}
\label{eqq10}
&      &   c_p\|u_1-u_2\|^p_{W^{m_1.p}(V)}+c_q\|v_1-v_2\|^q_{W^{m_2.q}(V)}\nonumber\\
& \geq &  \min\{c_p,c_q\}\left(\|u_1-u_2\|^q_{W^{m_1.p}(V)}+\|v_1-v_2\|^q_{W^{m_2.q}(V)}\right)\nonumber\\
& \geq & \frac{ \min\{c_p,c_q\}}{2^{q-1}}\|(u_1-u_2,v_1-v_2)\|^q\nonumber\\
& \geq &  \frac{\min\{c_p,c_q\}}{2^{p-1}}\|(u_1-u_2,v_1-v_2)\|^q.
\end{eqnarray}
Let
\begin{eqnarray}
\label{eqq5}
a_1(t)=\frac{\min\{c_p,c_q\}}{2^{p-1}}t^{q-1},\;\;\;\;t>2.
\end{eqnarray}
\par
(2) Assume that $\|u_1-u_2\|_{W^{m_1.p}(V)}\leq1$ and $\|v_1-v_2\|_{W^{m_2.q}(V)}\leq1$.  Then $\|(u_1-u_2,v_1-v_2)\| \leq 2$ and
\begin{eqnarray}
\label{eqq11}
&       &  c_p\|u_1-u_2\|^p_{W^{m_1.p}(V)}+c_q\|v_1-v_2\|^q_{W^{m_2.q}(V)}\nonumber\\
& \geq  & \min\{c_p,c_q\}\left(\|u_1-u_2\|^p_{W^{m_1.p}(V)}+\|v_1-v_2\|^p_{W^{m_2.q}(V)}\right)\nonumber\\
& \geq  & \frac{\min\{c_p,c_q\}}{2^{p-1}}\|(u_1-u_2,v_1-v_2)\|^p\nonumber\\
& \geq  & \begin{cases}
\frac{\min\{c_p,c_q\}}{2^{p-1}}\|(u_1-u_2,v_1-v_2)\|^p,\;\;\;\;\|(u_1-u_2,v_1-v_2)\|\leq1,\\
\frac{\min\{c_p,c_q\}}{2^{p-1}}\|(u_1-u_2,v_1-v_2)\|^q,\;\;\;\;1<\|(u_1-u_2,v_1-v_2)\|\leq2.
\end{cases}
\end{eqnarray}
Let
\begin{eqnarray}
\label{eqq6}
a_2(t)=
\begin{cases}
\frac{\min\{c_p,c_q\}}{2^{p-1}}t^{p-1},\;\;\;\;t\leq1,\\
\frac{\min\{c_p,c_q\}}{2^{p-1}}t^{q-1},\;\;\;\;1<t\leq2.
\end{cases}
\end{eqnarray}
\par
(3) Assume that $\|u_1-u_2\|_{W^{m_1.p}(V)}>1$ and $\|v_1-v_2\|_{W^{m_2.q}(V)} \leq 1$. Then $\|(u_1-u_2,v_1-v_2)\|^q > 1$ and
\begin{eqnarray}
\label{eqq12}
&      &  c_p\|u_1-u_2\|^p_{W^{m_1.p}(V)}+c_q\|v_1-v_2\|^q_{W^{m_2.q}(V)}\nonumber\\
& \geq &  \min\{c_p,c_q\}\left(\|u_1-u_2\|^q_{W^{m_1.p}(V)}+\|v_1-v_2\|^q_{W^{m_2.q}(V)}\right)\nonumber\\
& \geq &  \frac{\min\{c_p,c_q\}}{2^{q-1}}\|(u_1-u_2,v_1-v_2)\|^q\nonumber\\
& \geq &  \frac{\min\{c_p,c_q\}}{2^{p-1}}\|(u_1-u_2,v_1-v_2)\|^q.
\end{eqnarray}
Let
\begin{eqnarray}
\label{eqq7}
a_3(t)=\frac{\min\{c_p,c_q\}}{2^{p-1}}t^{q-1},\;\;\;\;t>1.
\end{eqnarray}
\par
(4) Assume that $\|u_1-u_2\|_{W^{m_1.p}(V)}\leq1$ and $\|v_1-v_2\|_{W^{m_2.q}(V)}>1$. Then $\|(u_1-u_2,v_1-v_2)\|^q>1$. Note that $q-p\le 0$. Thus, we have
\begin{eqnarray}
\label{eqq13}
&      &  c_p\|u_1-u_2\|^p_{W^{m_1.p}(V)}+c_q\|v_1-v_2\|^q_{W^{m_2.q}(V)}\nonumber\\
& \geq & \min\left\{c_p,c_q\|v_1-v_2\|^{q-p}_{W^{m_2.q}(V)}\right\}\left(\|u_1-u_2\|^p_{W^{m_1.p}(V)}+\|v_1-v_2\|^p_{W^{m_2.q}(V)}\right)\nonumber\\
& \geq & \min\left\{c_p,c_q(\|u_1-u_2\|_{W^{m_1.p}(V)}+\|v_1-v_2\|_{W^{m_2.q}(V)})^{q-p}\right\}
\left(\|u_1-u_2\|^p_{W^{m_1.p}(V)}+\|v_1-v_2\|^p_{W^{m_2.q}(V)}\right)\nonumber\\
& \geq & \min\left\{c_p,c_q\|(u_1-u_2,v_1-v_2)\|^{q-p}\right\}\frac{1}{2^{p-1}}\|(u_1-u_2,v_1-v_2)\|^p\nonumber\\
&  =   & \min\left\{\frac{c_p}{2^{p-1}}\|(u_1-u_2,v_1-v_2)\|^p,
\frac{c_q}{2^{p-1}}\|(u_1-u_2,v_1-v_2)\|^q\right\}\nonumber\\
&  \geq  & \frac{\min\{c_p,c_q\}}{2^{p-1}}\|(u_1-u_2,v_1-v_2)\|^q.
\end{eqnarray}
Let
\begin{eqnarray}
\label{eqq8}
a_4(t)=\frac{\min\{c_p,c_q\}}{2^{p-1}}t^{q-1},\;\;\;\;t>1.
\end{eqnarray}
Combining (\ref{eqq5}), (\ref{eqq6}), (\ref{eqq7}) and (\ref{eqq8}),  we define $a:\R^+\to \R^+$  by
\begin{eqnarray}
\label{eqq9}
a(t)=
\begin{cases}
\frac{\min\{c_p,c_q\}}{2^{p-1}}t^{p-1},\;\;\;\;t\leq1,\\
\frac{\min\{c_p,c_q\}}{2^{p-1}}t^{q-1},\;\;\;\;t>1.
\end{cases}
\end{eqnarray}
Then $a$ is continuous and strictly monotone increasing with $a(0)=0$ and $a(t)\to+\infty$ as $t\to +\infty$. Thus, by (\ref{eqq10}), (\ref{eqq11}), (\ref{eqq12}) and (\ref{eqq13}), (\ref{eqq4}) can be written as
\begin{eqnarray*}
\langle\Phi'(u_1,v_1)-\Phi'(u_2,v_2),(u_1-u_2,v_1-v_2)\rangle\\
\geq a(\|(u_1-u_2,v_1-v_2)\|_W)\|(u_1-u_2,v_1-v_2)\|_W.
\end{eqnarray*}
  So $\Phi'$ is uniformly monotone in $W$ if $\max\{p,q\}=p$. Similarly, if we let $\max\{p,q\}=q$, we can also obtain the same conclusion.
\par
Next, we show that $\Phi'$ is also hemicontinuous in $W$. Assume that $s\rightarrow s^*$ and $s,s^*\in[0,1]$. Note that
\begin{eqnarray*}
&        &  |\langle\Phi'((u_1,u_2)+s(v_1,v_2)),(w_1,w_2)\rangle-\langle\Phi'((u_1,u_2)+s^*(v_1,v_2)),(w_1,w_2)\rangle|\\
&  \leq  &    \|\Phi'((u_1,u_2)+s(v_1,v_2))-\Phi'((u_1,u_2)+s^*(v_1,v_2))\|_*\cdot\|(w_1,w_2)\|
\end{eqnarray*}
for all $(u_1,u_2),(v_1,v_2),(w_1,w_2)\in W$, where $\|\cdot\|_*$ denotes the norm of the dual space $W^*$. Then the continuity of $\Phi'$ implies that
$$
\langle\Phi'((u_1,u_2)+s(v_1,v_2)),(w_1,w_2)\rangle\rightarrow\langle\Phi'((u_1,u_2)+s^*(v_1,v_2))),(w_1,w_2)\rangle, \mbox{ as}\  s\rightarrow s^*
$$
for all $(u_1,u_2),(v_1,v_2),(w_1,w_2)\in W$. Hence $\Phi'$ is hemicontinuous in $W$.
\par
Moreover, for all $(u,v)\in W$, we have
\begin{eqnarray*}
\langle\Phi'(u,v),(u,v)\rangle
&=&\int_V\left[|\nabla^{m_1}u|^p+h_1(x)|u|^p\right]d\mu+\int_V\left[|\nabla^{m_2}v|^q+h_2(x)|v|^q\right]d\mu\\
&=&\|u\|^p_{W^{m_1,p}(V)}+\|v\|^q_{W^{m_2,q}(V)}.
\end{eqnarray*}
So $\Phi'$ is coercive in $W$. Thus by Theorem 26.A in \cite{E. Zeidler 1990}, we can obtain that $\Phi'$ admits a continuous inverse in $W$.
 \qed

\vskip2mm
\noindent
{\bf Lemma 3.5.} {\it $\Phi:W\to \R$ is sequentially weakly lower semi-continuous.}
\vskip0mm
\noindent
{\bf Proof.} Since $\Phi$ is continuously differentiable and $\Phi'$ is uniformly monotone (and then it is monotone), thus it follows from Proposition 25.20 in \cite{E. Zeidler 1990} that  $\Phi$ is sequentially weakly lower semi-continuous. \qed

\vskip2mm
\noindent
{\bf Lemma 3.6.} {\it $\Psi$ has compact derivative.}
\vskip0mm
\noindent
{\bf Proof.} Obviously, $\Psi$ is a $C^1$ functional on $W$. Assume that $\{(u_n,v_n)\}\subset W$ is bounded. Note that $W$ is of finite dimension. Then there exists a subsequence $\{(u_k,v_k)\}$ such that $(u_k,v_k)\rightarrow (u_0,v_0)$ for some $(u_0,v_0)\in W$. By the continuity of $\Psi'$, it is easy to obtain that
\begin{eqnarray*}
\|\Psi'(u_k,v_k)-\Psi'(u_0,v_0)\|_*\rightarrow 0, \;\mbox{as}\; n\rightarrow\infty.
\end{eqnarray*}
Hence, $\Psi'$ is compact in $W$.
\qed
\vskip2mm
\noindent
{\bf Proof of Theorem 1.1.} Obviously, by $(F_0)$ and $(F_1)$, $\Phi(0)=\Psi(0)=0$ and both $\Phi$ and $\Psi$ are continuously differentiable. Moreover, it is easy to see that $\Phi:W\to \R$ is coercive.
 Lemma 3.2-Lemma 3.6 imply that all of other conditions in Lemma 2.4 are satisfied. Hence, Lemma 2.4 implies that for each $\lambda\in\left( \Lambda_2^{-1},  \Lambda_1^{-1} \right)$,  the functional $\varphi$ has at least three distinct critical points that are solutions of  system (\ref{eq2}).
\qed

\vskip2mm
{\section{Proofs for the $(p,q)$-Laplacian system (\ref{p1})}}
  \setcounter{equation}{0}
Let $G=(V,E)$ be a locally finite graph. In order to investigate the $(p,q)$-Laplacian system (\ref{p1}), we work in the space $W_1:=W^{1,p}(V)\times W^{1,q}(V)$ with the norm  endowed with $\|(u,v)\|_1=\|u\|_{W^{1,p}(V)}+\|v\|_{W^{1,q}(V)}$ and then $(W_1,\|\cdot\|_1)$ is a  Banach space which is of infinite dimension. Different from the case of finite graph in section 2,  the continuous differentiability of variational functional for  (\ref{p1}) can not be obtained just by $(F_0)$. However, by using the condition  $(F_0)'$, the difficulty has been overcome in \cite{Yang 2023} so that we can  apply Lemma 2.4 to system (\ref{p1}).
\par
We consider the functional $\bar \varphi:W_1\to\R$ as
\begin{eqnarray}
\label{EQ1} \bar\varphi(u,v)=\frac{1}{p}\int_V(|\nabla u|^p+h_1(x)|u|^p)d\mu+\frac{1}{q}\int_V(|\nabla v|^q+h_2(x)|v|^q)d\mu-\lambda\int_V F(x,u,v)d\mu.
\end{eqnarray}
Then by Appendix A.2 in \cite{Yang 2023}, under the assumptions of Theorem 1.2, we have $\bar\varphi\in C^1(W_1,\mathbb{R})$, and
\begin{eqnarray}\label{EQ2}
\langle\bar\varphi'(u,v),(\phi_1,\phi_2)\rangle&=&\int_V\left[|\nabla u|^{p-2}\Gamma(u,\phi_1)+h_1(x)|u|^{p-2}u\phi_1-\lambda  F_u(x,u,v)\phi_1\right]d\mu\nonumber\\
&&+\int_V\left[(|\nabla v|^{q-2}\Gamma(v,\phi_2)+h_2(x)|v|^{q-2}v\phi_2-\lambda F_v(x,u,v)\phi_2\right]d\mu
\end{eqnarray}
for any $(u,v),(\phi_1,\phi_2)\in W_1$.  Define $\bar \Phi:W_1\rightarrow\mathbb{R}$ and $\bar \Psi:W_1\rightarrow\mathbb{R}$ by
  \begin{eqnarray*}
\bar\Phi(u,v)&=&\frac{1}{p}\int_V(|\nabla u|^p+h_1(x)|u|^p)d\mu+\frac{1}{q}\int_V(|\nabla v|^q+h_2(x)|v|^q)d\mu\\
&=&\frac{1}{p}\|u\|^p_{W^{1,p}(V)}+\frac{1}{q}\|v\|^q_{W^{1,q}(V)}
  \end{eqnarray*}
  and
\begin{eqnarray*}
\bar \Psi(u,v)=\int_V F(x,u,v)d\mu.
\end{eqnarray*}
Then $\bar\varphi(u,v)=\bar\Phi-\lambda\bar\Psi$. Moreover, it is easy to see that $(u,v)\in W_1$ is a critical point of $\bar\varphi$ if and only if
\begin{eqnarray*}
\int_V\left[|\nabla u|^{p-2}\Gamma(u,\phi_1)+h_1(x)|u|^{p-2}u\phi_1-\lambda F_u(x,u,v)\phi_1\right]d\mu=0
\end{eqnarray*}
 and
\begin{eqnarray*}
\int_V\left[|\nabla v|^{q-2}\Gamma(v,\phi_2)+h_2(x)|v|^{q-2}v\phi_2-\lambda F_v(x,u,v)\phi_2\right]d\mu=0
\end{eqnarray*}
for all $(\phi_1,\phi_2)\in W_1$.

 \vskip2mm
 \noindent
{\bf Lemma 4.1.}  {\it Assume that  $(M)$, $(H_1)$ and $(F_0)'$ hold. Then for any given $r>0$, the following inequality holds:
$$
\frac{\sup_{(u,v)\in\bar\Phi^{-1}(-\infty,r]}\bar\Psi(u,v)}{r}\le \frac{1}{r}\max_{|(s,t)|\le \frac{1}{h_0^{1/p}\mu_0^{1/p}} (pr)^{\frac{1}{p}}+ \frac{1}{h_0^{1/q}\mu_0^{1/q}} (qr)^{\frac{1}{q}}}a(|(s,t)|)\int_V b(x)d\mu.
$$ }
 \vskip0mm
 \noindent
{\bf Proof.} By $(F_0)'$, we have
\begin{eqnarray*}
        \bar \Psi(u,v)&=&\int_V{F(x,u(x),v(x))}d\mu\\
& \leq & \int_V a(|(u(x),v(x))|) b(x)d\mu\\
& \leq &  \max_{|(s,t)|\le\|u\|_\infty+\|v\|_\infty}a(|(s,t)|)\int_V b(x)d\mu
\end{eqnarray*}
for every $(u,v)\in W_1$. Furthermore,  for all $(u,v)\in W_1$ with $\bar\Phi(u,v)\leq r$, by Lemma 2.3, we get
\begin{eqnarray*}
\label{a1}
\|u\|_\infty \le \frac{1}{h_0^{1/p}\mu_0^{1/p}} \|u\|_{W^{1,p}(V)}\le \frac{1}{h_0^{1/p}\mu_0^{1/p}} (pr)^{\frac{1}{p}},\ \  \|v\|_\infty \le \frac{1}{h_0^{1/q}\mu_0^{1/q}} \|v\|_{W^{1,q}(V)}\le \frac{1}{h_0^{1/q}\mu_0^{1/q}} (qr)^{\frac{1}{q}}.
\end{eqnarray*}
Then
\begin{eqnarray}
\label{(8)}
\sup_{(u,v)\in\bar\Phi^{-1}(-\infty,r]}\bar\Psi(u,v)\leq   \max_{|(s,t)|\le \frac{1}{h_0^{1/p}\mu_0^{1/p}} (pr)^{\frac{1}{p}}+ \frac{1}{h_0^{1/q}\mu_0^{1/q}} (qr)^{\frac{1}{q}}}a(|(s,t)|)\int_V b(x)d\mu.
\end{eqnarray}
Furthermore, the proof is completed by multiplying $\frac{1}{r}$ in both side of (\ref{(8)}). \qed

\vskip2mm
 \noindent
{\bf Lemma 4.2.} {\it Assume that $(M)$, $(H_1)$, $(F_0)'$ and $(F_3)'$ hold. Then there exists $ (u_{\delta_1},v_{\delta_2})\in W_1$  such that}
$$
 \frac{\sup_{(u,v)\in\bar\Phi^{-1}(-\infty,\gamma_1^{p}+\gamma_2^{q}]}\bar\Psi(u,v)}{\gamma_1^{p}+\gamma_2^{q}}
 < \frac{\bar\Psi(u_{\delta_1},v_{\delta_2})}{\bar\Phi(u_{\delta_1},v_{\delta_2})}.
$$
\vskip0mm
 \noindent
{\bf Proof.}   Let
\begin{eqnarray*}
u_{\delta_{1}}(x)=\begin{cases}
                  \delta_1,& x=x_0\\
                  0,&x\not= x_0
                  \end{cases},
\quad v_{\delta_{2}}(x)=\begin{cases}
                  \delta_2,& x=x_0\\
                  0,&x\not= x_0,
                  \end{cases}
\end{eqnarray*}
where $\delta_i, i=1,2$ are defined in $(F_3)'$. Then a simple calculation implies that
$$
|\nabla u_{\delta_1}|(x)=\begin{cases}
           \sqrt{\frac{deg(x_0)}{2\mu(x_0)}}\delta_1, & x=x_0,\\
           \sqrt{\frac{w_{x_0y}}{2\mu(y)}}\delta_1, & x=y \mbox{ with } y\thicksim x_0,\\
           0, &\mbox{otherwise},
           \end{cases}
$$
and
$$
|\nabla v_{\delta_2}|(x)=\begin{cases}
           \sqrt{\frac{deg(x_0)}{2\mu(x_0)}}\delta_2, & x=x_0,\\
           \sqrt{\frac{w_{x_0y}}{2\mu(y)}}\delta_2, & x=y \mbox{ with } y\thicksim x_0,\\
           0, &\mbox{otherwise}.
           \end{cases}
$$
Then
\begin{eqnarray}\label{pp1}
&     &     \int_V{(|\nabla u_{\delta_1}|^p+h_1(x)|u_{\delta_1}|^p)}d\mu\nonumber\\
&  =  &    \sum_{x\in V}(|\nabla u_{\delta_1}|^p+h_1(x)|u_{\delta_1}|^p)\mu(x)\nonumber\\
&  =  &    (|\nabla u_{\delta_1}|^p(x_0)+h_1(x_0)|u_{\delta_1}|^p(x_0))\mu(x_0)+\sum_{y\thicksim x_0}(|\nabla u_{\delta_1}|^p(y)+ h_1(y)|u_{\delta_1}|^p(y))\mu(y)\nonumber\\
&  =  & \left(\frac{deg(x_0)}{2\mu(x_0)}\right)^{\frac{p}{2}}\delta_1^p\mu(x_0)+h_1(x_0)\delta_1^p\mu(x_0)+\delta_1^p\sum_{y\thicksim x_0}\left(\frac{w_{x_0y}}{2\mu(y)}\right)^{\frac{p}{2}}\mu(y)\nonumber\\
&  =  & \delta_1^p M_1 \quad (\mbox{defined in } (F_3)'),
\end{eqnarray}
and similarly,
\begin{eqnarray}\label{pp2}
&     &     \int_V{(|\nabla v_{\delta_2}|^q+h_2(x)|v_{\delta_2}|^q)}d\mu\nonumber\\
&  =  & \left(\frac{deg(x_0)}{2\mu(x_0)}\right)^{\frac{q}{2}}\delta_2^q\mu(x_0)+h_2(x_0)\delta_2^q\mu(x_0)+\delta_2^q\sum_{y\thicksim x_0}\left(\frac{w_{x_0y}}{2\mu(y)}\right)^{\frac{q}{2}}\mu(y)\nonumber\\
&  =  & \delta_2^qM_2 \quad (\mbox{defined in } (F_3)').
\end{eqnarray}
Note that $\{y|y\thicksim x_0\}$ is a finite set. Then (\ref{pp1}) and  (\ref{pp2}) imply that  $(u_{\delta_1},v_{\delta_2})\in W_1$.
Moreover,
\begin{eqnarray}
\label{(5)}
    \bar\Phi(u_{\delta_1},v_{\delta_2})
 & = &  \frac{1}{p}\int_V{(|\nabla u_{\delta_1}|^p+h_1(x)|u_{\delta_1}|^p)}d\mu+\frac{1}{q}\int_V{(|\nabla v_{\delta_2}|^q+h_2(x)|v_{\delta_2}|^q)}d\mu\nonumber\\
 & = &   \frac{\delta_1^pM_1}{p}+\frac{\delta_2^qM_2}{q}.
\end{eqnarray}
Note that $\delta_i>\gamma_i\kappa_i$, $i=1,2$, where
$$
\kappa_1= \left(\frac{M_1}{p}\right)^{-\frac{1}{p}}, \ \ \kappa_2=\left(\frac{M_2}{q}\right)^{-\frac{1}{q}}.
$$
Then (\ref{(5)}) implies that $ \bar\Phi(u_{\delta_1},v_{\delta_2})\ge \gamma_1^p+\gamma_2^q$.
 Moreover, $(F_1)'$ implies that
\begin{eqnarray}
\label{(6)}
        \bar\Psi(u_{\delta_1},v_{\delta_2})&=&\int_VF(x,u_{\delta_1},v_{\delta_2})d\mu\nonumber\\
&  =   & F(x_0,\delta_1,\delta_2)+ \int_{V/\{x_0\}}F(x,0,0)d\mu\nonumber\\
&  =   & F(x_0,\delta_1,\delta_2).
\end{eqnarray}
Hence, by (\ref{(5)}) and (\ref{(6)}), we get
\begin{eqnarray*}
      \frac{\bar\Psi(u_{\delta_1},v_{\delta_2})}{\bar\Phi(u_{\delta_1},v_{\delta_2})}
=  \frac{F(x_0,\delta_1,\delta_2)}{ \frac{\delta_1^pM_1}{p}+\frac{\delta_2^qM_2}{q}}.
\end{eqnarray*}
In view of Lemma 4.1 and $(F_3)'$, we get
\begin{eqnarray*}
&      &  \frac{\sup_{(u,v)\in\bar\Phi^{-1}(-\infty,\gamma_1^{p}+\gamma_2^{q}]}\bar\Psi(u,v)}{\gamma_1^{p}+\gamma_2^{q}}\\
& \leq & \frac{1}{\gamma_1^{p}+\gamma_2^{q}}\max_{|(s,t)|\le \frac{1}{h_0^{1/p}\mu_0^{1/p}} (p\gamma_1^{p}+p\gamma_2^{q})^{\frac{1}{p}}+ \frac{1}{h_0^{1/q}\mu_0^{1/q}} (q\gamma_1^{p}+q\gamma_2^{q})^{\frac{1}{q}}}a(|(s,t)|)\int_V b(x)d\mu \\
& < & \frac{F(x_0,\delta_1,\delta_2)}{ \frac{\delta_1^pM_1}{p}+\frac{\delta_2^qM_2}{q}}\\
& = & \frac{\bar\Psi(u_{\delta_1},v_{\delta_2})}{\bar\Phi(u_{\delta_1},v_{\delta_2})}.
\end{eqnarray*}
The proof is complete. \qed

\vskip2mm
\noindent
{\bf Lemma 4.3.} {\it Assume that $(H_1)$ and  $(F_2)'$ hold. Then for each $\lambda\in \left(0,+\infty\right)$, the functional $\bar\Phi-\lambda\bar\Psi$ is coercive.}
\vskip0mm
\noindent
{\bf Proof.} By $(\bar F_2)'$, we have
\begin{eqnarray*}
       \varphi(u,v)
&  =   & \frac{1}{p}\|u\|_{W^{1,p}(V)}^p+\frac{1}{q}\|v\|_{W^{1,q}(V)}^q-\lambda\int_V{F(x,u,v)}d\mu\\
& \ge  & \frac{1}{p}\|u\|_{W^{1,p}(V)}^p+\frac{1}{q}\|v\|_{W^{1,q}(V)}^q-\lambda \|f_1\|_\infty\|u\|_{\infty}^{\alpha}-\lambda \|f_2\|_\infty\|v\|_{\infty}^\beta-\lambda\int_Vg(x)d\mu\\
& \ge  & \frac{1}{p}\|u\|_{W^{1,p}(V)}^p+\frac{1}{q}\|v\|_{W^{1,q}(V)}^q-\lambda \|f_1\|_\infty h_0^{-\frac{\alpha}{p}}\mu_0^{-\frac{\alpha}{p}}\|u\|_{W^{1,p}(V)}^\alpha\\
&      &  -\lambda \|f_2\|_\infty h_0^{-\frac{\beta}{q}} \mu_0^{-\frac{\beta}{q}}\|v\|_{W^{1,q}(V)}^\beta -\lambda\int_Vg(x)d\mu.
\end{eqnarray*}
Note that $\alpha\in[0,p)$ and $\beta\in[0,q)$. Therefore $\bar\varphi(u,v)$ is a coercive functional for every $\lambda\in \left(0,+\infty\right)$.
\qed

\vskip2mm
\noindent
{\bf Lemma 4.4.} {\it The G$\hat{a}$teaux derivative of $\bar\Phi$ admits a continuous inverse on $W_1^*$, where  $W_1^*$ is the dual space of $W_1$.}
\vskip0mm
\noindent
{\bf Proof.} In the proof of Lemma 3.4, we only need to take $m_i=1, i=1,2$ and let $G=(V,E)$ be a locally finite graph. Then the proof is the essentially same as that in Lemma 3.4.
 \qed

\vskip2mm
\noindent
{\bf Lemma 4.5.} {\it $\bar\Phi:W_1\to \R$ is sequentially weakly lower semi-continuous.}
\vskip0mm
\noindent
{\bf Proof.} The proof is the essentially same as that in Lemma 3.5 with  taking $m_i=1, i=1,2$ and letting $G=(V,E)$ be a locally finite graph.\qed

\vskip2mm
\noindent
{\bf Lemma 4.6.} {\it $\bar\Psi$ has compact derivative.}
\vskip0mm
\noindent
{\bf Proof.} Obviously, $\bar\Psi$ is a $C^1$ functional on $W_1$. Assume that $\{(u_n,v_n)\}\subset W_1$ is bounded. Then by Lemma 2.3, there exists a positive constant $M>0$ such that $\|u_k\|_\infty\le M$ and $\|v_k\|_\infty\le M$ and there exists a subsequence $\{(u_k,v_k)\}$ such that $(u_k,v_k)\rightharpoonup (u_0,v_0)$ for some $(u_0,v_0)\in W_1$.  In particular,
$$\lim_{k\rightarrow\infty}\int_Vu_k\varphi d\mu=\int_Vu_0\varphi d\mu,\forall\varphi\in C_c(V),$$
which implies that
\begin{eqnarray}\label{eq14}
\lim_{k\rightarrow\infty}u_k(x)=u_0(x)\;\; \mbox{for any fixed}\;x\in V
\end{eqnarray}
by taking
$$
\varphi(y)=
\begin{cases}
1,& y=x\\
0, & y\not=x.
\end{cases}
$$
Similarly, we have
\begin{eqnarray*}\label{m1}
\lim_{k\rightarrow\infty}v_k(x)=v_0(x)\;\; \mbox{for any fixed}\;x\in V.
\end{eqnarray*}
Note that
\begin{eqnarray*}
&     &\|\bar\Psi'(u_k,v_k)-\bar\Psi'(u_0,v_0)\|_*\\
&  =  &\sup_{\|(\phi_1,\phi_2)\|=1}|\langle\bar\Psi'(u_k,v_k)-\bar\Psi'(u_0,v_0),(\phi_1,\phi_2)\rangle|\\
&  =  &\sup_{\|(\phi_1,\phi_2)\|=1}|\int_V[F_{u_k}(x,u_k,v_k)-F_{u_0}(x,u_0,v_0)]\phi_1d\mu
+\int_V[F_{v_k}(x,u_k,v_k)-F_{v_0}(x,u_0,v_0)]\phi_2d\mu|\\
&\leq &\sup_{\|(\phi_1,\phi_2)\|=1}|\int_V[F_{u_k}(x,u_k,v_k)-F_{u_0}(x,u_0,v_0)]\phi_1d\mu|
+|\int_V[F_{v_k}(x,u_k,v_k)-F_{v_0}(x,u_0,v_0)]\phi_2d\mu|.
\end{eqnarray*}
By $(F_0)'$, we have
\begin{eqnarray*}
&     &|F_{u_k}(x,u_k,v_k)-F_{u_0}(x,u_0,v_0)|\\
&\leq & [a(|(u_k,v_k)|)+a(|(u_0,v_0)|)]b(x)\\
&\leq & \left[\max_{|(s,t)|\leq \|u_0\|_\infty+\|v_0\|_\infty}a(|(s,t)|)+\max_{|(s,t)|\leq 2M}a(|(s,t)|)\right]b(x)\\
&:= & l(x).
\end{eqnarray*}
 Note that $b\in L^1(V)$. Hence, $l(x)\in L^1(V)$ and so $\int_V|F_{u_k}(x,u_k,v_k)-F_{u_0}(x,u_0,v_0)|d\mu$ is uniformly convergent. Thus, by (\ref{eq14}), (\ref{m1}) and the continuity of $F_u$, we  have
\begin{eqnarray*}
&     &\int_V[F_{u_k}(x,u_k,v_k)-F_{u_0}(x,u_0,v_0)]\phi_1d\mu\\
&\leq &\left(\int_V|F_{u_k}(x,u_k,v_k)-F_{u_0}(x,u_0,v_0)|d\mu\right)\|\phi_1\|_\infty\\
&\leq &\left(\int_V|F_{u_k}(x,u_k,v_k)-F_{u_0}(x,u_0,v_0)|d\mu\right)h_0^{-\frac{1}{p}}l_0^{-\frac{1}{p}}\\
&\rightarrow & 0\;\;\mbox{as}\; k\rightarrow\infty.
\end{eqnarray*}
Similarly, we also have
\begin{eqnarray*}
\int_V[F_{v_k}(x,u_k,v_k)-F_{v_0}(x,u_0,v_0)]\phi_2d\mu
\rightarrow  0\;\;\mbox{as}\; k\rightarrow\infty.
\end{eqnarray*}
So,
\begin{eqnarray*}
&     &\|\bar\Psi'(u_k,v_k)-\bar\Psi'(u_0,v_0)\|_*\rightarrow 0\;\;\mbox{as}\;k\rightarrow\infty.
\end{eqnarray*}
Hence, $\bar\Psi'$ is compact in $W_1$.
\qed
\vskip2mm
\noindent
{\bf Proof of Theorem 1.2.} Obviously, by $(F_0)'$ and $(F_1)'$, $\bar\Phi(0)=\bar\Psi(0)=0$ and both $\bar\Phi$ and $\bar\Psi$ are continuously differentiable. Moreover, it is easy to see that $\bar\Phi:W_1\to \R$ is coercive.
 Lemma 4.2-Lemma 4.6 imply that all of other conditions in Lemma 2.4 are satisfied.  Hence, Lemma 2.4 implies that for each $\lambda\in\left( \Theta_2^{-1},  \Theta_1^{-1} \right)$,  the functional $ \bar\varphi$ has at least  three distinct critical points that are solutions of  system (\ref{p1}).
\qed

\vskip2mm
\noindent
{\bf Remark 4.1.} On the locally finite graph, we do not consider the more general poly-Laplacian system. That is because it is difficult to obtain the continuous differentiability of the variational functional $\varphi$ when $m_i>1, i=1,2$, which is caused by the special definition of $\mathcal{L}_{m,p}$.

\vskip2mm
{\section{The results for the scalar equation}}
 \setcounter{equation}{0}
\vskip2mm
\par
By using the similar arguments of Theorem 1.1, we can also obtain a similar result for the following scalar equation on finite graph $G=(V,E)$:
\begin{eqnarray}
\label{eqq21}
  \pounds_{m,p}u+h(x)|u|^{p-2}u=\lambda f(x,u),\;\;\;\;\hfill x\in V,
\end{eqnarray}
where $m\geq 1$ is an integer, $h:V\rightarrow\mathbb{R}$, $p>1$, $\lambda>0$ and  $f:V\times\mathbb{R}\rightarrow\mathbb{R}$.

\vskip2mm
\noindent
{\bf Theorem 5.1.} {\it Let $G=(V,E)$ be a finite graph and $F(x,s)=\int_0^s f(x,\tau)d\tau$ for all $x\in V$. Assume that the following conditions hold: \\
$(h)$\; $h(x)>0$ for all $x\in V$;\\
$(f_0)$\; $F(x,s)$ is continuously differentiable in $s\in \R$ for all $x\in V$;\\
$(f_1)$\; $\int_VF(x,0)d\mu=0;$\\
$(f_2)$ \; there exists a constant $\alpha \in [0,p)$ and functions $g_1,g_2:V\to \R$ such that
$$
F(x,s)\le g_1(x)|s|^\alpha+g_2(x)
$$
 for all $(x,s)\in V\times\mathbb{R}$;\\
$(f_3)$ \; there are  positive constants $\gamma$ and $\delta$ with $\delta>\gamma\kappa$, such that
$$
 \Lambda_1':=\frac{1}{\gamma^{p}}\max_{x\in V,|s| \le{\frac{(p{\gamma}^{p})^{\frac{1}{p}}}{h_{\min}^{1/p}\mu_{\min}^{1/p}} }}F(x,s)|V|\\
< \frac{\inf_{x\in V}F(x,\delta)|V|}{ \frac{\delta^p}{p}\int_Vh(x)d\mu}:=\Lambda_2',\\
$$
where $|V|=\sum_{x\in V}\mu(x)$ and  $\kappa=\left(\frac{1}{p}\int_Vh(x)d\mu\right)^{-\frac{1}{p}}$.\\
  Then  for each parameter $\lambda$ belonging to
$\left( \Lambda_2'^{-1},  \Lambda_1'^{-1} \right),$
 system (\ref{eqq21}) has at least three distinct solutions.}

\vskip2mm
\par
By using the similar arguments of Theorem 1.2, we can also obtain a similar result for the following scalar equation on a locally finite graph $G=(V,E)$:
\begin{eqnarray}
\label{eqq2}
   -\Delta_pu+h(x)|u|^{p-2}u=\lambda f(x,u),\;\;\;\;\hfill x\in V,
\end{eqnarray}
where $p\ge 2$, $h:V\rightarrow\mathbb{R}$, $\lambda>0$  and $f:V\times\mathbb{R}\rightarrow\mathbb{R}$.

\vskip2mm
\noindent
{\bf Theorem 5.2.} {\it Let $G=(V,E)$ be a locally finite graph and $F(x,s)=\int_0^s f(x,\tau)d\tau$ for all $x\in V$. Assume that $(M)$ and the following conditions hold: \\
$(h)'$\; there exists a constant $h_0>0$ such that $h(x)\geq h_0>0$ for all $x\in V$;\\
$(f_0)'$\; $F(x,s)$ is continuously differentiable in $s\in \R$ for all $x\in V$,  and there exists a  function $a\in C(\R^+,\R^+)$ and a function $b:V\to \R^+$ with $b\in L^1(V)$ such that
$$
|f(x,s)|\le a(|s|) b(x),  |F(x,s)|\le a(|s|) b(x),
$$
for all $x\in V$ and all $s\in \R$;\\
$(f_1)'$\; $\int_V F(x,0)d\mu=0$ and there exists a $x_0\in V$ such that $F(x_0,0)=0$;\\
$(f_2)'$ \; there exists a constants $\alpha \in [0,p)$ and functions $g_i:V\to \R$, $i=1,2$, with $g_1\in L^\infty(V)$ and $g_2\in L^1(V)$, such that
$$
F(x,s)\le g_1(x)|s|^\alpha+g_2(x)
$$
 for all $(x,s)\in V\times\mathbb{R}$;\\
$(f_3)'$ \; there are positive constants $\gamma$ and $\delta$ with $\delta>\gamma\kappa$, such that
$$
 \Theta_1':=\frac{1}{\gamma^{p}}\max_{|s|\le \frac{1}{h_0^{1/q}\mu_0^{1/q}} (p\gamma^{p})^{\frac{1}{p}}}a(|s|)\int_V b(x)d\mu \\
< \frac{F(x,\delta)d\mu}{ \frac{\delta^p}{p}M}:=\Theta_2',\\
$$
where $\kappa=\left(\frac{1}{p}\int_Vh(x)d\mu\right)^{-p}$ and
$$ M =  \left(\frac{deg(x_0)}{2\mu(x_0)}\right)^{\frac{p}{2}}\mu(x_0)+h(x_0)\mu(x_0)+\sum_{y\thicksim x_0}\left(\frac{w_{x_0y}}{2\mu(y)}\right)^{\frac{p}{2}}\mu(y).$$
  Then  for each parameter $\lambda$ belonging to
$\left( \Theta_2'^{-1},  \Theta_1'^{-1} \right),$
 system (\ref{eqq2}) has at least three  solutions.}

\vskip2mm
\noindent
{\bf Remark 5.1.} In Theorem 5.1 and Theorem 5.2, all of  three solutions are nontrivial solutions if we furthermore assume that
$f(x,0)\not=0$ for some $x\in V$.

\vskip2mm
{\section{Examples}}
 \setcounter{equation}{0}
\vskip2mm
 \noindent
 In this section, we present two examples as applications of Theorem 1.1 and Theorem 5.2.
\vskip2mm
 \noindent
 {\bf Examples 6.1.} Let $p=2, q=3$ and $m=2$ in (\ref{eq2}). Consider the following system:
 \begin{eqnarray}
\label{6.1}
 \begin{cases}
  \pounds_{2,2}u+h_1(x)u=\lambda F_u(x,u,v),\;\;\;\;\hfill x\in V,\\
  \pounds_{2,3}v+h_2(x)v=\lambda F_v(x,u,v),\;\;\;\;\hfill x\in V,\\
   \end{cases}
\end{eqnarray}
where $G=(V,E)$  is a finite graph of $9$ vertexes, that is,  $V=\{x_1,x_2,\cdots ,x_9\}$, the measure $\mu(x_i)=1,i=1,2,\cdots,9,h_i:V\rightarrow\mathbb{R}^+$ with $h_i\equiv9,i=1,2,$ for all $x\in V$ and put
$$
\omega_1=\frac{(p\gamma_1^{p}+p\gamma_2^{q})^{\frac{1}{p}}}{h_{1,\min}^{1/p}\mu_{\min}^{1/p}},
\omega_2=\frac{(q\gamma_1^{p}+q\gamma_2^{q})^{\frac{1}{q}}}{h_{2,\min}^{1/q}\mu_{\min}^{1/q}}.
$$
Let $\lambda>0$ and $F:V\times\R\times\R\to\R$ is  defined by
\begin{eqnarray*}
\frac{\partial{F(x,s,t)}}{{\partial s}}=
\begin{cases}
\omega_1-|s|,    &0\leq|s|\leq \omega_1,\\
|s|^3-\omega_1^3,&\omega_1<|s|<4\omega_1,\\
(4\omega_1)^{r_1}|s|^{3-r_1}-\omega_1^3,&4\omega_1\leq|s|,
\end{cases}
\end{eqnarray*}
and
\begin{eqnarray*}
\frac{\partial{F(x,s,t)}}{{\partial t}}=
\begin{cases}
\omega_2-|t|,&0\leq|t|\leq\omega_2,\\
|t|^4-\omega_2^4,&\omega_2<|t|<5\omega_2,\\
(5\omega_2)^{r_2}|t|^{4-r_2}-\omega_2^4,&5\omega_2\leq|t|.
\end{cases}
\end{eqnarray*}
Then
{\small \begin{eqnarray*}
F(x,s,t)=
\begin{cases}
\frac{1}{2}\omega_1^2+\frac{1}{4}(4\omega_1)^4+\frac{3}{4}\omega_1^4+\frac{1}{4-r_1}(4\omega_1)^{r_1}|s|^{4-r_1}-\omega_1^3|s|-\frac{1}{4-r_1}(4\omega_1)^4\\
+\frac{1}{2}\omega_2^2+\frac{1}{5}(5\omega_2)^5+\frac{4}{5}\omega_2^5+\frac{1}{5-r_2}(5\omega_2)^{r_2}|t|^{5-r_2}-\omega_2^4|t|-\frac{1}{5-r_2}(5\omega_2)^5,
&4\omega_1\leq|s|,5\omega_2\leq|t|,\\
\omega_1|s|-\frac{1}{2}|s|^2+\frac{1}{2}\omega_2^2+\frac{1}{5}|t|^5-\omega_2^4|t|+\frac{4}{5}\omega_2^5,&0\leq|s|\leq \omega_1,\omega_2<|t|<5\omega_2,\\
\omega_1|s|-\frac{1}{2}|s|^2+\frac{1}{2}\omega_2^2+\frac{1}{5}(5\omega_2)^5+\frac{4}{5}\omega_2^5+\\
\frac{1}{5-r_2}(5\omega_2)^{r_2}|t|^{5-r_2}-\omega_2^4|t|-\frac{1}{5-r_2}(5\omega_2)^5,
&0\leq|s|\leq \omega_1,5\omega_2\leq|t|,\\
\frac{1}{2}\omega_1^2+\frac{1}{4}|s|^4-\omega_1^3|s|+\frac{3}{4}\omega_1^4+\omega_2|t|-\frac{1}{2}|t|^2,&\omega_1<|s|<4\omega_1,0\leq|t|\leq \omega_2,\\
\frac{1}{2}\omega_1^2+\frac{1}{4}|s|^4-\omega_1^3|s|+\frac{3}{4}\omega_1^4+\frac{1}{2}\omega_2^2+\frac{1}{5}|t|^5-\omega_2^4|t|+\frac{4}{5}\omega_2^5,
&\omega_1<|s|<4\omega_1,\omega_2<|t|<5\omega_2,\\
\frac{1}{2}\omega_1^2+\frac{1}{4}|s|^4-\omega_1^3|s|+\frac{3}{4}\omega_1^4
+\frac{1}{2}\omega_2^2+\frac{1}{5}(5\omega_2)^5+\frac{4}{5}\omega_2^5+\frac{1}{5-r_2}(5\omega_2)^{r_2}|t|^{5-r_2}\\
-\omega_2^4|t|-\frac{1}{5-r_2}(5\omega_2)^5,
&\omega_1<|s|<4\omega_1,5\omega_2\leq|t|,\\
\frac{1}{2}\omega_1^2+\frac{1}{4}(4\omega_1)^4+\frac{3}{4}\omega_1^4+\frac{1}{4-r_1}(4\omega_1)^{r_1}|s|^{4-r_1}-\omega_1^3|s|-\frac{1}{4-r_1}(4\omega_1)^4\\
+\omega_2|t|-\frac{1}{2}|t|^2,&4\omega_1\leq|s|,0\leq|t|\leq \omega_2,\\
\frac{1}{2}\omega_1^2+\frac{1}{4}(4\omega_1)^4+\frac{3}{4}\omega_1^4+\frac{1}{4-r_1}(4\omega_1)^{r_1}|s|^{4-r_1}-\omega_1^3|s|-\frac{1}{4-r_1}(4\omega_1)^4\\
+\frac{1}{2}\omega_2^2+\frac{1}{5}|t|^5-\omega_2^4|t|+\frac{4}{5}\omega_2^5,&4\omega_1\leq|s|,\omega_2<|t|<5\omega_2,\\
\omega_1|s|-\frac{1}{2}|s|^2+\omega_2|t|-\frac{1}{2}|t|^2,&0\leq|s|\leq \omega_1,0\leq|t|\leq \omega_2,
\end{cases}
\end{eqnarray*}}
where $(r_1,r_2)\in(1,2]\times(1,3]$. Next we verify that $h_1,h_2$ and $F$ satisfy the conditions in Theorem 1.1.
 \par
 $\bullet$ Obviously, $h_i$ satisfy $(H)$, $i=1,2$, and $F$ satisfies $(F_0)$ and $(F_1)$.
 \par
 $\bullet$ Let
 $$f_1(x)\equiv \frac{(4\omega_1)^{r_1}}{4-r_1},f_2(x)\equiv \frac{(5\omega_2)^{r_2}}{5-r_2}, $$
 and
 $$
 g(x)\equiv \frac{1}{2}\omega_1^2+\frac{1}{4}(4\omega_1)^4+\frac{3}{4}\omega_1^4+\frac{1}{2}\omega_2^2+\frac{1}{5}(5\omega_2)^5+\frac{4}{5}\omega_2^5, \mbox{ for all }x\in V.
 $$
 Then
 $$F(x,s,t)\leq  f_1(x)|s|^\alpha+f_2(x)|t|^\beta+g(x),$$
 where $\alpha\in [0,p),\beta\in[0,q).$ Hence,  $F$ satisfies $(F_2)$.
 \par
 $\bullet$ Let
 $$\delta_1=4\omega_1=4\times15^\frac{1}{2},\delta_2=5\omega_2=5\times(\frac{45}{2})^\frac{1}{3},\gamma_1=(\frac{81}{2})^\frac{1}{2},\gamma_2=(\frac{81}{3})^\frac{1}{3}.$$
Then $\delta_1>\gamma_1\kappa_1=1$, $\delta_2>\gamma_2\kappa_2=1$,
 \begin{eqnarray*}
 \Lambda_2^{-1}
 &  =  &\frac{ \frac{\delta_1^p}{p}\int_Vh_1(x)d\mu+\frac{\delta_2^q}{q}\int_Vh_2(x)d\mu}{\inf_{x\in V}F(x,\delta_1,\delta_2)|V|}\\
 &  =  &\frac{\frac{81}{2}\cdot (4\cdot15^\frac{1}{2})^2+\frac{81}{3}\cdot(5\cdot(\frac{45}{2})^\frac{1}{3})^3}{\left(\frac{1}{2}\cdot15^\frac{2}{2}+\frac{1}{4}(4\cdot15^\frac{1}{2})^4-4\cdot15^\frac{4}{2}+\frac{3}{4}\cdot15^\frac{4}{2}+
\frac{1}{2}\cdot(\frac{45}{2})^\frac{2}{3}+\frac{1}{5}(5\cdot(\frac{45}{2})^\frac{1}{3})^5-5\cdot(\frac{45}{2})^\frac{5}{3}+\frac{4}{5}\cdot(\frac{45}{2})^\frac{5}{3}\right)\times 9}\\
&\approx & 0.07614,
 \end{eqnarray*}
and
 \begin{eqnarray*}
 \Lambda_1^{-1}
  &  =  &\frac{\gamma_1^{p}+\gamma_2^{q}} {\max_{x\in V,|s|\le{\frac{(p\gamma_1^{p}+p\gamma_2^{q})^{\frac{1}{p}}}{h_{1,\min}^{1/p}\mu_{\min}^{1/p}} }, |t|\le \frac{ (q\gamma_1^{p}+q\gamma_2^{q})^{\frac{1}{q}}}{h_{2,\min}^{1/q}\mu_{\min}^{1/q}}}F(x,s,t)|V|}\\
  &  =  &\frac{\frac{81}{2}+\frac{81}{3}}{\left(\frac{1}{2}\cdot(15)^\frac{2}{2}+\frac{1}{2}\cdot(\frac{45}{2})^\frac{2}{3}\right)\times 9}\\
  &\approx & 0.65303>\Lambda_2^{-1}.
 \end{eqnarray*}
 Hence,  $(F_3)$ holds.  Thus, by Theorem 1.1, for each $\lambda\in ( \Lambda_2^{-1}, \Lambda_1^{-1})\approx(0.07614,0.65303)$, system (\ref{6.1}) has at least three distinct solutions.

 \vskip2mm
\vskip2mm
 \noindent
 {\bf Examples 6.2.} Let $p=3$  in (\ref{eqq2}). Consider the following scalar equation on locally finite graph $G=(V,E)$:
\begin{eqnarray}
\label{6.2}
  -\Delta_3u+h(x)u=\lambda f(x,u),\;\;\;\;\hfill x\in V,
\end{eqnarray}
where the measure $\mu(x)\equiv 1$ and $h(x)\equiv 4$ for all $x\in V$. For some fixed $x_0\in V$, there are $4$ edges $x_0y\in E$ with $w_{x_0y}=2$. Put
$$
\omega=\frac{(p\gamma^{p})^{\frac{1}{p}}}{h_0^{1/p}\mu_0^{1/p}}.
$$
Let $\lambda>0$ and $F:V\times\R\to \R$ is  defined by
\begin{eqnarray*}
f(x_0,s)=\frac{\partial F(x_0,s)}{\partial s}=
\begin{cases}
\omega-|s|,&|s|\leq \omega,\\
|s|^5-\omega^5,&\omega<|s|\leq 6\omega,\\
(6\omega)^rs^{5-r}-\omega^5,&6\omega<|s|,
\end{cases}
\end{eqnarray*}
and
$$f(x,s)=0,\;\;\;\mbox{for all}\;x\in V/\{x_0\}.$$
Then
\begin{eqnarray*}
F(x_0,s)=
\begin{cases}
\omega|s|-\frac{1}{2}|s|^2,&|s|\leq \omega,\\
\frac{1}{2}\omega^2+\frac{1}{6}|s|^6-\omega^5|s|+\frac{5}{6}\omega^6,&\omega<|s|\leq 6\omega,\\
\frac{1}{2}\omega^2+\frac{6^6+5}{6}\omega^6-\frac{1}{6-r}(6\omega)^6+\frac{1}{6-r}(6\omega)^r|s|^{6-r}-\omega^5|s|,&6\omega<|s|
\end{cases}
\end{eqnarray*}
and
$$F(x,s)=0,\;\;\;\mbox{for all}\;x\in V/\{x_0\},$$
 where $r\in(3,5]$ and $\lambda>0$. Next we verify that $h$ and $F$ satisfy the conditions in Theorem 5.2.
 \par
 $\bullet$ Obviously, $h$ satisfy $(h)'$.
 \par
 $\bullet$ Let
  \begin{eqnarray*}
 g_1(x)=
 \begin{cases}
 \frac{1}{6-r}(6\omega)^r,&x=x_0,\\
 0,&x\neq x_0,
 \end{cases}
 \end{eqnarray*}
   \begin{eqnarray*}
 g_2(x)=
 \begin{cases}
\frac{1}{2}\omega^2+\frac{6^6+5}{6}\omega^6,&x=x_0,\\
 0,&x\neq x_0,
 \end{cases}
 \end{eqnarray*}
   \begin{eqnarray*}
 a(|s|)=
 \begin{cases}
\omega|s|-\frac{1}{2}|s|^2+1,&|s|\leq \omega,\\
\frac{1}{2}\omega^2+\frac{1}{6}|s|^6-\omega^5|s|+\frac{5}{6}\omega^6+1,&\omega<|s|\leq 6\omega,\\
\frac{1}{2}\omega^2+\frac{6^6+5}{6}\omega^6-\frac{1}{6-r}(6\omega)^6+\frac{1}{6-r}(6\omega)^r|s|^{6-r}-\omega^5|s|+1,&6\omega<|s|,
\end{cases}
 \end{eqnarray*}
and
 \begin{eqnarray*}
 b(x)=
 \begin{cases}
 1,&x=x_0,\\
 0,&x\neq x_0.
 \end{cases}
 \end{eqnarray*}
 Then
 $$\|g_1\|_\infty=\frac{1}{6-r}(6\omega)^r,\|g_2\|_{L^1(V)}=\frac{1}{2}\omega^2+\frac{6^6+5}{6}\omega^6,$$
 and
 $$F(x,s)\leq g_1(x)|s|^\alpha+g_2(x).$$
Moreover,
 $$f(x,s)\leq a(|s|)b(x),F(x,s)\leq a(|s|)b(x),$$
 for all $x\in V $ and all $s\in \R$. Hence, $F$ satisfies $(f_0)',(f_1)'$ and $(f_2)'$.
 \par
 $\bullet$ Let
 $$\delta=6\omega=6\cdot4^\frac{1}{3},\gamma=(\frac{16}{3})^\frac{1}{3}.$$
Then $\delta>\gamma\kappa=1$,
 \begin{eqnarray*}
\Theta_2^{-1}
 &  =  &\frac{ \frac{\delta^p}{p}M}{\inf_{x\in V}F(x,\delta)}\\
 &  =  &\frac{(6\cdot4^\frac{1}{3})^3\times \frac{1}{3}\times16}{\frac{1}{2}\cdot4^\frac{2}{3}+6^5\cdot4^2-6\cdot4^2+\frac{5}{6}\cdot4^2}\\
 &\approx & 0.0371
 \end{eqnarray*}
and
 \begin{eqnarray*}
\Theta_1^{-1}
  &  =  &\frac{\gamma^{p}} {\max_{x\in V,|s|\le{\frac{(p\gamma^{p})^{\frac{1}{p}}}{h_{1,\min}^{1/p}\mu_{\min}^{1/p}} }}a(|s|)\int_Vb(x)d\mu}\\
  &  =  &\frac{\frac{16}{3}}{\frac{1}{2}\cdot4^\frac{2}{3}+1}\\
  &\approx &2.36 >\Theta_2^{-1}.
 \end{eqnarray*}
 Hence, $F$ satisfy $(f_3)'$.  Thus, by Theorem 5.2, for each $\lambda\in ( \Theta_1^{-2}, \Theta_1^{-1})\approx(0.0371,2.36)$, equation (\ref{6.2}) has at least three distinct  solutions.
 \vskip2mm
 \noindent
 {\bf Acknowledgments}\\
This project is supported by Yunnan Fundamental Research Projects (grant No: 202301AT070465) and  Xingdian Talent Support Program for Young Talents of Yunnan Province.

\vskip2mm
 \noindent
 {\bf Authors' contributions}\\
Yan Pang and Xingyong Zhang contribute the main manuscript equally.

\vskip2mm
\renewcommand\refname{References}
{}
\end{document}